\documentclass[12pt,leqno]{article}
\usepackage{amssymb,amsmath,amsfonts,amsbsy, xspace,
latexsym, amscd, graphicx}

\newtheorem{theorem}{Theorem}[section]
\newtheorem{lemma}{Lemma}[section]

\newtheorem{definition}{Definition}[section]

\newtheorem{simpler notation}{Simpler
Notation}[section]
\newtheorem{remark}{Remark}[section]
\newtheorem{example}{Example}[section]

\newcommand{\ifff}{{\em\, i\,f\,{f}\quad}}

\newcommand{\bra}{\ensuremath{\large\langle}}

\newcommand{\ket}{\ensuremath{\large\rangle}}

\newcommand{\supp}{\ensuremath{{\rm{supp}}}}

\begin{document}

\pagenumbering{roman}
\title{Non-Standard Analysis, Multiplication of
Schwartz Distributions, and Delta-Like Solution of Hopf's
Equation\\} \date{} \maketitle
\begin{center}
A Thesis Presented to the Faculty of\\
California Polytechnic State University\\
                        San Luis Obispo\\[1.5in]

In Partial Fulfillment of the Requirements for the Degree\\
Master of Science in Mathematics \\[3cm]
by\\
Guy Berger\\
September 2005
\end{center}

\newpage
\begin{center}
AUTHORIZATION FOR REPRODUCTION OF MASTER'S THESIS\\[2cm]
\end{center}
\begin{flushleft}
I grant permission for the reproduction of this thesis in its
entirety or any of its parts, without further authorization from
me.\\[5cm]
\end{flushleft}
$\overline{\text{Signature \quad\quad\quad\quad\quad
}}$\\\\\\\\
$\overline{\text{Date \quad\quad\quad\quad\quad\quad\quad  }}$

\newpage
\begin{center}
APPROVAL PAGE \\[2cm]
\end{center}
\begin{flushleft}
TITLE: Non-Standard Analysis, Multiplication of Schwartz
Distributions, and Delta-Like Solution of Hopf's Equation\\
AUTHOR: Guy Berger\\
DATE SUBMITTED: September 16, 2005\\[3cm]
\end{flushleft}
\begin{tabbing}
\= $\overline{\text{Adviser\quad\quad\quad\quad\quad\quad\quad\quad}}$\hspace{3cm}\= $\overline{\text{Signature\quad\quad\quad\quad\quad}}$\\[3cm]
\> $\overline{\text{Committee Member\quad\quad\,\;\;}}$ \> $\overline{\text{Signature\quad\quad\quad\quad\quad}}$\\[3cm]
\> $\overline{\text{Committee Member\quad\quad\,\;\;}}$ \> $\overline{\text{Signature\quad\quad\quad\quad\quad}}$\\[3cm]

\end{tabbing}

\begin{center}
Acknowledgements\\[2cm]
\end{center}
\begin{flushleft}
I would like to thank the faculty and staff of the Cal Poly
Mathematics Department, who have made my experience over the past
two years a positive one.\\
I am especially grateful to my supervisor, Dr. Todor D. Todorov,
for his help, guidance, and support throughout the past year.
\end{flushleft}
\newpage
\begin{abstract}
\begin{center}
\textbf{Non-Standard Analysis, Multiplication of Schwartz
Distributions, and Delta-Like Solution of Hopf's Equation}\\
Guy Berger \end{center} We construct an algebra of generalized
functions $^*\mathcal{E}(\mathbb{R}^d)$. We also construct an
embedding of the space of  Schwartz distributions
$\mathcal{D}^\prime(\mathbb{R}^d)$ into
$^*\mathcal{E}(\mathbb{R}^d)$ and thus present a solution of the
problem of multiplication of Schwartz distributions which improves
J.F. Colombeau's solution. As an application we prove the
existence of a weak delta-like solution in
${^*\mathcal{E}(\mathbb{R}^d)}$ of the Hopf equation. This
solution does not have a counterpart in the classical theory of
partial differential equations. Our result improves a similar
result by M. Radyna obtained in the framework of perturbation
theory.
\end{abstract}

\newpage
\tableofcontents

\newpage

\pagenumbering{arabic}
\section{Introduction}\label{S: Intro}

\quad\quad     In what follows
$\mathcal{E}(\mathbb{R}^d)=\mathcal{C}^\infty(\mathbb{R}^d)$
denotes the class of $\mathcal{C}^\infty$-functions on
$\mathbb{R}^d$. Also
$\mathcal{D}(\mathbb{R}^d)=\mathcal{C}_0^\infty(\mathbb{R}^d)$
denotes the class of test functions on $\mathbb{R}^d$ and
$\mathcal{D}^\prime(\mathbb{R}^d)$ stands for the space of
Schwartz distributions (Schwartz generalized functions) on
$\mathbb{R}^d$ (H. Bremermann~ \cite{hBremermann65}).

    The algebra of generalized functions
${^*\mathcal{E}(\mathbb{R}^d)}$ is a particular non-standard
extension of the class $\mathcal{E}(\mathbb{R}^d)$. The field of
the scalars $^*\mathbb{C}$ of the algebra
${^*\mathcal{E}(\mathbb{R}^d)}$ is a particular non-standard
extension of the field of complex numbers $\mathbb{C}$ and the
field of the real scalars $^*\mathbb{R}$ is a non-standard
extension of  $\mathbb{R}$. That means that both $^*\mathbb{C}$
and $^*\mathbb{R}$ are non-Archimedean fields containing non-zero
infinitesimals, i.e. generalized numbers $h$ such that $0<|h|<1/n$
for all $n\in\mathbb{N}$. Since the involvement of non-Archimedean
fields in applied mathematics is somewhat unusual, we start with a
summary  of the relevant definitions and results in the theory of
ordered fields and non-Archimedean fields (Section 2).

    In Sections~\ref{S: Filters and
Ultrafilters}-\ref{S: Ultrafilter on D(Rd)}  we present the basic
facts of the theory of free filters and ultrafilters (C. C. Chang
and H. J. Keisler~\cite{CKeis}). We construct a particular
ultrafilter on the space of test functions
$\mathcal{D}(\mathbb{R}^d)$ which is important for the embedding
of Schwartz distributions in the algebra
${^*\mathcal{E}(\mathbb{R}^d)}$.

    In Sections 5-6 we present the construction of the
fields of the complex and real non-standard numbers $^*\mathbb{C}$
and $^*\mathbb{R}$. In Section~\ref{S: Saturation Principle in *C}
we prove the {\bf Saturation Principle in $^*\mathbb{C}$} which
plays a role in non-standard analysis similar to the role of the
completeness of $\mathbb{R}$ and $\mathbb{C}$ in usual (standard)
analysis. These sections might be viewed as an introduction to
non-standard analysis (A. Robinson~\cite{aRob66}). We should note
that our exposition of non-standard analysis does not require any
background in mathematical logic or model theory.

    The construction of the algebra
${^*\mathcal{E}(\mathbb{R}^d)}$ is presented in Section~\ref{S:
Non-Standard Smooth Functions}; in short,
${^*\mathcal{E}(\mathbb{R}^d)}$ is a {\bf differential associative
commutative algebra of generalized functions} similar to (but much
larger than) the class
$\mathcal{E}(\mathbb{R}^d)=\mathcal{C}^\infty(\mathbb{R}^d)$.
 In
Section~\ref{S: Internal Sets and Saturation Principle in *E} we
state the Saturation Principle for
${^*\mathcal{E}(\mathbb{R}^d)}$, playing the role of the
completeness property.

    In Section~\ref{S: Embedding of Schwartz
Distributions
     in *E} we construct the chain of embeddings
$\mathcal{E}(\mathbb{R}^d)\subset\mathcal{D}^\prime(\mathbb{R}^d)
\subset{^*\mathcal{E}(\mathbb{R}^d)}$. These embeddings presents a
{\bf solution of the problem of multiplication of Schwartz
distributions} similar to but different from Colombeau's solution
of the same problem (J.F. Colombeau~\cite{jfCol84}). The problem
of multiplication of Schwartz distributions has an interesting and
dramatic history. Soon after the distribution theory was invented
by L. Schwartz, he proved that the space of distributions
$\mathcal{D}^\prime(\mathbb{R}^d)$ cannot be supplied with an
associative and commutative product that reproduces the usual
product in the spaces $\mathcal{C}^k(\mathbb{R}^d), k=0, 1,
2,\dots$. This negative result, known as {\em Schwartz
Impossibilities Result }(L. Schwartz~\cite{l.Schwartz54}), was the
reason this problem was considered for a long time as unsolvable.
In the late 1980's Jean F. Colombeau offered a solution of the
problem of multiplication of distributions by constructing an
algebra of generalized functions $\mathcal{G}(\mathbb{R}^d)$ with
the chain of algebraic embeddings
$\mathcal{E}(\mathbb{R}^d)\subset\mathcal{D}^\prime(\mathbb{R}^d)
\subset\mathcal{G}(\mathbb{R}^d)$ thus avoiding Schwartz
Impossibilities Result (since $k=\infty$). One (slightly
disturbing) feature of Colombeau's solution is that the set of
scalars $\overline{\mathbb{C}}$ of the algebra
$\mathcal{G}(\mathbb{R}^d)$ is a ring with zero divisors, not a
field as any set of scalars should be. In this respect our
solution of the problem of multiplication of Schwartz
distributions presents an important {\bf  improvement of
Colombeau's theory}: the set of scalars $^*\mathbb{C}$ of the
algebra $^*\mathcal{E}(\mathbb{R}^d)$ is an {\bf algebraically
complete $c^+$-saturated field} (Section~\ref{S: Saturation
Principle in *C}). As a consequence, the set of the real scalars
$^*\mathbb{R}$ is a {\bf real closed Cantor complete field}. We
should notice that the fact that $^*\mathcal{E}(\mathbb{R}^d)$ is
a differential algebra (not merely a linear space) is important
for our goals in applied mathematics, in particular, for studying
generalized solutions of non-linear partial differential equations
such as shock-wave and delta-like solutions. Notice that these are
the solutions after the formation of the shock in many
conservation law type equations.

    In Section~\ref{S: Generalized Delta-Like Solution
of Hopf Equation} we  prove the existence of a weak delta-like
solution of the Hopf equation $u_t(x, t)+u(x, t)u_x(x, t) = 0$ in
the framework of ${^*\mathcal{E}(\mathbb{R}^d)}$. This solution
has counterparts neither in the spaces of classical functions such
as $\mathcal{C}^k(\mathbb{R}^d),\; k=1, 2, \dots, \infty$, nor in
the spaces of Schwartz distributions such as
$\mathcal{D}^\prime(\mathbb{R}^d)$. Our result improves a similar
result by M. Radyna~\cite{mRadynaICGF2000} obtained in the spirit
of perturbation theory.
\newpage


\section{Ordered Fields}\label{S: Ordered Fields}
We will begin by defining ordered fields and giving examples of
some well known and some lesser known orderings on the
(non-archimedean) fields of rational functions and Laurent
series.\\
    \begin{definition}\label{D:Ordering}
    Let $\mathbb{K}$ be a field (ring).
    $\mathbb{K}$ is called \textbf{orderable}
    if there exists a nonempty subset
    $\mathbb{K}_{+}\subset\mathbb{K}$
    such that
    \begin{itemize}
        \item[(1)] $0\notin\mathbb{K}_{+}$
        \item[(2)] $x,y\in\mathbb{K}_{+}
\Longrightarrow
        x+y,xy\in\mathbb{K}_{+}$
        \item[(3)] For every non-zero
$x\in\mathbb{K}$,  either $x\in\mathbb{K}_{+}$ or
$-x\in\mathbb{K}_{+}$
    \end{itemize}
\end{definition}
$\mathbb{K}_{+}$ generates an order relation $<_{\mathbb{K}}$ on
$\mathbb{K}$ as follows: $x<_{\mathbb{K}}y$ \,\ifff\,
$y-x\in\mathbb{K}_{+}.$ \\\\
 $(\mathbb{K},<_{\mathbb{K}})$ is called a
\emph{totally ordered
field} or simply an \emph{ordered field}.\\\\
We will write $<$ instead of $<_{\mathbb{K}}$ when it is clear
from context which field's order relation we are
referring to.\\\\
In addition to $(\mathbb{R},<)$ (where $<$ is the usual order on
$\mathbb{R}$) there are many (more interesting) examples of
ordered fields. But first let us make a short detour:
\begin{example}\label{E:C Unorderable}
$\mathbb{C}$, the set of complex numbers, is not orderable.
    \end{example}
    \emph{Proof} Suppose there exists a subset
$\mathbb{C}_{+}$
    satisfying Definition \ref{D:Ordering}. Then
consider:\\
    \textbf{Case 1} Suppose $i\in\mathbb{C}_{+}$. Then
$i\cdot
    i=-1\in\mathbb{C}_{+}$, implying $(-1)\cdot
   ( -1)=1\in\mathbb{C}_{+}$. This is impossible since
    $-1+1=0\notin\mathbb{C}_{+}$.\\
    \textbf{Case 2} Suppose $i\notin\mathbb{C}_{+}$.
Then
    $-i\in\mathbb{C}_{+}$, implying $(-i)\cdot
    (-i)=-1\in\mathbb{C}_{+}$, leading to the same
contradiction as
    in \textbf{Case 1}.  $\blacktriangle$ \\\\
    The previous example can be generalized as
follows:
    \begin{theorem}\label{T:Characterization}
    A field $\mathbb{K}$ is orderable \ifff it is
formally real. This means that for every
    $n\in\mathbb{N}$ and every $x_{k}\in\mathbb{K}$\\
    \[
        \sum_{k=1}^{n} x_{k}^{2}=0\quad implies\quad
x_{k}=0\quad
        for\quad all\quad k.\]
    \end{theorem}
    For details on the subject of formally real fields
and the proof of this theorem,
    see (Van Der Waerden~\cite{VanDerWaerden}, Chapter
11).

    \begin{definition}\label{D:Ordered field
homomorphism}
    Let $\mathbb{K}$ and $\mathbb{L}$ be ordered
fields. If
    $\varphi:\;\mathbb{K}\longrightarrow\mathbb{L}$ is
a field
    homomorphism that preserves order, i.e.
$x<_{\mathbb{K}}y$
    implies $\varphi(x)<_{\mathbb{L}}\varphi(y)$, then
$\varphi$
    is said to be an \textbf{ordered field
homomorphism}.\\
    \textbf{Ordered field isomorphisms} and
\textbf{ordered field
    embeddings} are defined similarly.

    \end{definition}

    \begin{remark}\label{R:Q embedding}
        There exists an ordered field embedding from
$\mathbb{Q}$
        into any ordered field $\mathbb{K}$. We call
it the
        \textbf{canonical embedding} of $\mathbb{Q}$
into
        $\mathbb{K}$ and it is defined by:
$\sigma(0)=0$,
        $\sigma(n)=n\cdot1$ and
$\sigma(-n)=-\sigma(n)$ for
        $n\in\mathbb{N}$, and
$\sigma(p/q)=\sigma(p)/\sigma(q)$
        for $p,\,q\in\mathbb{Z},\,q\ne0$.
    \end{remark}
    From now on, if $x\in\mathbb{N},\mathbb{Z},$ or
$\mathbb{Q}$,
    we will refer to $x$ and $\sigma(x)\in\mathbb{K}$
interchangeably.

    \begin{example}\label{E:Rational Functions}
    Let $(\mathbb{R},<)$ be the field of real numbers
with the usual order,
    and let $\mathbb{R}(x)$ be the set
    of rational functions in the variable $x$ with
coefficients in
    $\mathbb{R}$.  Note that we may think of
$\mathbb{R}$ as a subfield
    of $\mathbb{R}(x)$, as represented by the constant
functions. Then define
    $\mathbb{R}(x)_{+}=\{R(x):R(x)\in\mathbb{R}(x)$
and there
    exists $x_{0}\in\mathbb{R}$ such that $R(x)>0$
whenever
    $x>x_{0}\}.$\\
    The ordered field generated by
$\mathbb{R}(x)_{+}$, which we will refer to
    simply as $\mathbb{R}(x)$, has some surprising
properties.
    Namely:
        \begin{itemize}
            \item[(i)] $\mathbb{R}(x)$ contains
\textbf{infinitely large
            elements} like $f(x)=x$. This means that
            $f(x)>n$ for all $n\in\mathbb{N}$ (let
$x_{0}=n$).
            \item[(ii)] $\mathbb{R}(x)$ contains
\textbf{positive infinitesimals} like $g(x)=\frac{1}{x}$. This
means that
            $0<g(x)<\frac{1}{n}$ for all
$n\in\mathbb{N}$ (let $x_{0}=n$).

        \end{itemize}
        \end{example}

   \begin{remark}\label{R:Inherited order}
        Let $\mathbb{L}$ be an ordered integral domain
(an ordered
        ring without zero divisors) and $\mathbb{K}$
be the field
        of fractions of \, $\mathbb{L}$.
        Define
\[
        \mathbb{K}_{+}=\{\frac{x}{y}:x,y\in\mathbb{L}_{+}\;
        or\; -x,-y\in \mathbb{L}_{+}\}
        \]
        The order generated by $\mathbb{K}_{+}$ is the
only one
        which extends the order in $\mathbb{L}$. It is
said to be
        the order \textbf{inherited} from
$\mathbb{L}$.
    \end{remark}

    \begin{example}\label{E:Rational Functions
revisited}
    With Remark~(\ref{R:Inherited order}) in mind, we
may revisit
    Example~(\ref{E:Rational Functions}).\\
    If \; $\mathbb{R}[x]$ is
    the ring of polynomials over $\mathbb{R}$, we may
define
    $\mathbb{R}[x]_{+}=\{P(x)\in \mathbb{R}[x]$:
    $lead(P)>0\}$, where $lead(P)$ is the leading
    coefficient of $P(x)$.\\
    The order generated on $\mathbb{R}[x]$ by
$\mathbb{R}[x]_{+}$
    can be extended to $\mathbb{R}(x)$ since
$\mathbb{R}(x)$ is
    the field of fractions of $\mathbb{R}[x]$. That
is, we may
    redefine
$\mathbb{R}(x)_{+}=\{\frac{P(x)}{Q(x)}:P(x),Q(x)\in\mathbb{R}[x]
    \; and \; lead(P),lead(Q)>0 \; or \;
    lead(P),lead(Q)<0\}$\\
    This definition is equivalent to that given
previously and the orders
    generated by the two are in fact one and the same.
Now we may plainly see
    that $f(x)=x$ is indeed an infinitely large
element since
    $lead(x-n)=1>0$ for all $n\in\mathbb{N}$. \\
    Similarly, $g(x)=\frac{1}{x}$ is a positive
infinitesimal
    because $\frac{1}{n}-\frac{1}{x}=\frac{x-n}{nx}$
and
    $lead(x-n),lead(nx)>0$ for all $n\in\mathbb{N}$.

    \end{example}

    \begin{example}\label{E:Laurent Series}
    The set
    \[
\mathbb{R}(x^{\mathbb{Z}})=\{\sum_{n=m}^{\infty}\frac{a_{n}}{x^{n}}:\quad
    a_{n}\in\mathbb{R}\quad, \quad m\in\mathbb{Z},
\quad and \quad a_{m}\ne 0\}
    \]
     of \textbf{Laurent series} with
    coefficients in $\mathbb{R}$ is a field under
normal polynomial
addition and multiplication.\\
    We may define an order on
    $\mathbb{R}(x^\mathbb{Z})$ by
\[
\mathbb{R}_{+}(x^{\mathbb{Z}})=\{\sum_{n=m}^{\infty}\frac{a_{n}}{x^{n}}\in\mathbb{R}(x^{\mathbb{Z}}):\quad
        a_{m}>0\}
        \]
    Here, an element such as
    \[
\frac{1}{x}+\frac{1}{x^{2}}+\frac{1}{x^{3}}+\ldots
        \]
    is infinitesimal.
    \end{example}

\begin{example}\label{E:Laurent Series 2}
    The field of (formal) Laurent Series may also be
defined as
    follows:
    \[
\mathbb{R}(x^{\mathbb{Z}})=\{\sum_{n=m}^{\infty}a_{n}{x^{n}}:\quad
    a_{n}\in\mathbb{R},\quad m\in\mathbb{Z},\quad
and\quad a_{m}\ne0\}
    \]
    If we now let
    \[
\mathbb{R}_{+}(x^{\mathbb{Z}})=\{\sum_{n=m}^{\infty}
    a_{n}{x^{n}}\in\mathbb{R}(x^{\mathbb{Z}}):\quad
        a_{m}>0\}
        \]
then even a series that is divergent for all x, such as
\[
    x + 2x^{2} + 6x^{3} + \ldots + n!x^{n} +
    \ldots
\]
is an infinitesimal in this field.
\end{example}

    \section{Filters and Ultrafilters}\label{S:
Filters and Ultrafilters} In this section we define and give
examples of filters on an arbitrary infinite set. Having done so,
we will prove the existence of an ultrafilter using the axiom
of choice.\\
    \begin{definition}\label{D: filter}
    Let $I$ be an infinite set and let
$\mathcal{F}\subset\mathcal{P}(I),\;\mathcal{F}\ne\varnothing$. If
    $\mathcal{F}$ satisfies:
    \begin{itemize}
        \item[(F1)] If $A\in\mathcal{F}$ and $A\subset
B\subset I$,
        then $B\in\mathcal{F}$.
        \item[(F2)] $A,\,B\in\mathcal{F}$ implies
        $A\cap B\in\mathcal{F}$
        \item[(F3)] $\varnothing\notin\mathcal{F}$

    \end{itemize}
    then $\mathcal{F}$ is a \textbf{filter} on $I$.
    If it is also true that
    \begin{itemize}
        \item[(F4)]
$\bigcap_{A\in\mathcal{F}}A=\varnothing$
    \end{itemize}
    then $\mathcal{F}$ is called a \textbf{free
filter} on $I$.
    A filter $\mathcal{F}$ is called countably
incomplete if:
    \begin{itemize}
        \item[(F5)] There exists a sequence of
decreasing sets
        $I=I_{0}\supset I_{1} \supset I_{2} \supset
\ldots$ in
        $\mathcal{F}$ such that
        $\bigcap_{n=0}^{\infty}I_{n}=\varnothing$.
    \end{itemize}
    \end{definition}
If $\mathcal{F}$ is a filter on $I$, it follows immediately from
the definition that:

\begin{itemize}
    \item[(i)] $I\in\mathcal{F}$
    \item[(ii)] $\mathcal{F}$ is closed under finite
    intersections.
    \item[(iii)] If $A\in\mathcal{F}$ then $I\setminus
    A\notin\mathcal{F}$
    \item[(iv)] If $\mathcal{F}$ is countably
incomplete, then
    $\mathcal{F}$ is free.
\end{itemize}

\begin{definition}\label{D: Ultrafilter}
A filter $\mathcal{U}$ on a set $I$ is called an
\textbf{ultrafilter} if for every filter $\mathcal{F}$ on $I$,
$\mathcal{U}$ is not a proper subset of $\mathcal{F}$. That is,
there is no filter $\mathcal{F}$ on $I$ that properly contains
$\mathcal{U}$.
\end{definition}

\begin{theorem}\label{T: Ultrafilter characterization}
Let $\mathcal{U}$ be a filter on $I$. Then $\mathcal{U}$ is an
ultrafilter on $I$ \ifff for every $A\subset I$ either
$A\in\mathcal{U}$ or $I\setminus A\in\mathcal{U}$.
\end{theorem}
\emph{Proof} Suppose $\mathcal{U}$ is an ultrafilter on $I$ and
$A,I\setminus A \notin\mathcal{U}$. Let
$\widehat{\mathcal{U}}=\{X\; : \; A\cup X\in\mathcal{U}\}$. It is
not hard to check that $\widehat{\mathcal{U}}$ is a filter that
properly contains $\mathcal{U}$, since $I\setminus
A\in\widehat{\mathcal{U}}$. Thus $\mathcal{U}$ cannot be an
ultrafilter. To prove the other direction, suppose, to the
contrary, that $\mathcal{U}$ is not an ultrafilter. Then there
exists a filter $\mathcal{V}$ which is a proper extension of
$\mathcal{U}$. Let $A\in\mathcal{V}\setminus\mathcal{U}$. Then,
since $\mathcal{U}\subset\mathcal{V}$, we have that $A\cap
X\neq\varnothing$ for all $X\in\mathcal{U}$. But $I\setminus
A\in\mathcal{U}$ by assumption, so $A\cap \left(I\setminus
A\right)\neq\varnothing$, a contradiction.
$\blacktriangle$.\\

\begin{example}\label{E: Fixed filter}
Let $I=\mathbb{N}$, and fix $n\in\mathbb{N}$. Then
$\mathcal{U}=\{X\;:\;X\subset\mathbb{N},\; n\in X\}$ is an
ultrafilter on $\mathbb{N}$. However, $\mathcal{U}$ is clearly not
free.

\end{example}

\begin{example}\label{E: Frechet filter}
The filter $\mathcal{F}_{r}(\mathbb{N})$ consisting of all
cofinite sets of natural numbers is called the Fr\'{e}chet filter
on $\mathbb{N}$. It is free and countably incomplete since
$\bigcap_{n=1}^{\infty}\left(\mathbb{N}\setminus\{n\}\right)=\varnothing$
However, since neither the set of even numbers nor the set of odd
numbers is in $\mathcal{F}_{r}(\mathbb{N})$, by Theorem~(\ref{T:
Ultrafilter characterization}), it is not an ultrafilter.

\end{example}

\begin{theorem}\label{T: Ultrafilter existence}
Let $I$ be an infinite set and $\mathcal{F}$ a free filter on $I$.
Then there exists a free ultrafilter $\mathcal{U}$ on $I$ such
that $\mathcal{F}\subseteq\mathcal{U}$.
\end{theorem}
\emph{Proof} Let $\widehat{\mathcal{F}}$ be the set of all filters
on $I$ that contain $\mathcal{F}$.
$\widehat{\mathcal{F}}\neq\varnothing$ since
$\mathcal{F}\in\widehat{\mathcal{F}}$. Let $\widehat{\mathcal{F}}$
be ordered by set inclusion, and consider a linearly ordered
subset $\mathcal{M}\subseteq\widehat{\mathcal{F}}$. Define
$\widehat{\mathcal{M}}=\bigcup_{M\in\mathcal{M}}M$. Note that if
$A\in\widehat{\mathcal{M}}$ then $A\in M$ for some
$M\in\mathcal{M}$. Thus, if $A\subset B\subset I$, it follows that
$B\in M$, implying $B\in\widehat{\mathcal{M}}$. Also, if $A,B\in
\widehat{\mathcal{M}}$, then we must have that $A\in M_{1}$ and
$B\in M_{2}$ for some $M_{1},M_{2}\in \mathcal{M}$. Since
$\mathcal{M}$ is linearly ordered, we may assume without loss of
generality that $M_{1}\subset M_{2}$. Thus $A,B\in M_{2}$. Hence
$A\cap B\in M_{2}$, implying $A\cap B\in\widehat{\mathcal{M}}$.
Finally, we have that $\varnothing\notin\widehat{\mathcal{M}}$
because otherwise $\varnothing$ would be an element of some filter
$M\in\mathcal{M}$, which is impossible. We have just shown that
$\widehat{\mathcal{M}}$ is itself a filter. But since the choice
of $\mathcal{M}$ was arbitrary, we can conclude that
\textbf{every} linearly ordered subset of $\widehat{\mathcal{F}}$
has an upper bound in $\widehat{\mathcal{F}}$. Thus, by Zorn's
Lemma, $\widehat{\mathcal{F}}$ has a maximal element,
$\mathcal{U}$, which is an ultrafilter on $I$ containing
$\mathcal{F}$. Also, since $\mathcal{F}\subset\mathcal{U}$ and
$\mathcal{F}$ is free, we have that
$\bigcap_{A\in\mathcal{U}}A\subseteq\bigcap_{A\in\mathcal{F}}A=\varnothing$.
Therefore $\mathcal{U}$ is free. $\blacktriangle$.\\

    \section{Ultrafilter on
$\mathcal{D}(\mathbb{R}^{d})$}\label{S: Ultrafilter on D(Rd)} Here
we define an ultrafilter on $\mathcal{D}(\mathbb{R}^d)$, the set
of test functions, in order to construct ordered, non-archimedean
fields of non-standard real and complex numbers, $^*\mathbb{R}$
and $^*\mathbb{C}$,
respectively.\\
\begin{definition}\label{D: Basic Sets}
Let $\mathcal{D}(\mathbb{R}^{d})$ be the set of test functions on
$\mathbb{R}^{d}$. That is,
$\mathcal{D}(\mathbb{R}^{d})=\mathcal{C}^{\infty}_{0}(\mathbb{R}^{d})$.
For every $n\in\mathbb{N}$, define the \textbf{basic set}
$\mathcal{B}_{n}$ by
\begin{align*}
\mathcal{B}_{n}=\{&\varphi\in\mathcal{D}(\mathbb{R}^{d})\,:\
\\
                      & \varphi \text{ is real-valued
and symmetric},
\\
                      & \varphi(x)=0 \text{ for all }
x\in \mathbb{R}^{d}, \Vert x\Vert\geq
                      1/n,\\
                      & \int\varphi=1\\
                      & \int x^{\alpha}\varphi=0
\text{ for all } \alpha\in\mathbb{N}^{d}_{0},
                       1\leq\vert\alpha\vert\leq n,\\
                      &
1\leq\int\vert\varphi\vert<1+\frac{1}{n}\}
\end{align*}
and $\mathcal{B}_{0}=\mathcal{D}(\mathbb{R}^{d})$.

\end{definition}

\begin{theorem}\label{T: Properties of Basic Sets}
\begin{itemize}
    \item[(i)] $\mathcal{B}_{n}\neq\varnothing$ for
all n.
    \item[(ii)] $\mathcal{B}_{0}\supset\mathcal{B}_{1}
    \supset\mathcal{B}_{2}\supset\ldots$
    \item[(iii)]
$\bigcap_{n}\mathcal{B}_{n}=\varnothing$

\end{itemize}
\end{theorem}
\emph{Proof} For the proof of (i), see (Oberguggenberger and
Todorov~\cite{OberTod98}). (ii) follows from Definition~(\ref{D:
Basic Sets}). For (iii), suppose there were a function $\varphi$
such that $\varphi\in\bigcap_{n}\mathcal{B}_{n}$ for all $n$. Then
consider $\widehat{\varphi}(\xi)=\int \varphi(x)e^{i\xi x}dx$, the
Fourier transform of $\varphi$. Since
$\varphi\in\mathcal{D}(\mathbb{R}^{d})$, $\widehat{\varphi}$ is
entire (\emph{Bremermann}~\cite{hBremermann65} Lemma 8.11, p.85).
Therefore, we can write
$\widehat{\varphi}(\xi)=\sum_{\alpha\in\mathbb{N}_{0}^{d}}
\frac{(\partial^{\alpha}\widehat{\varphi})(0)}{\alpha!}\xi^{\alpha}$.
But $0=i^{\vert\alpha\vert}\int
x^{\alpha}\varphi(x)dx=\left.i^{\vert\alpha\vert}\int
x^{\alpha}\varphi(x)e^{i\xi x}dx\right|_{\xi=0} =
(\partial^{\alpha}\widehat{\varphi})(0)$ for all $\alpha\neq 0$.
It follows that $\widehat{\varphi}$ is constant. However, by the
same lemma as before, we also have that
$\lim_{\vert\xi\vert\to\infty}\widehat{\varphi}(\xi)=0$. Thus
$\widehat{\varphi}=0$, implying $\varphi(x)=0$. This contradicts
the property that $\int \varphi = 1$. Hence
$\bigcap_{n}\mathcal{B}_{n}=\varnothing$.
$\blacktriangle$.\\

\begin{definition}\label{D: Basic filter}
Define the \textbf{basic filter} $\mathcal{F}_{\mathcal{B}}$ on
$\mathcal{D}(\mathbb{R}^{d})$ by
\[
\mathcal{F}_{\mathcal{B}}=\{\Phi\subseteq\mathcal{D}(\mathbb{R}^{d})\;
    : \; \mathcal{B}_{n}\subseteq\Phi \text{ for some
}
    n\in\mathbb{N}\}.
    \]
\end{definition}
Since each $\mathcal{B}_{n}$ is itself an element of
$\mathcal{F}_{\mathcal{B}}$, it follows from Theorem~(\ref{T:
Properties of Basic Sets}) that $\mathcal{F}_{\mathcal{B}}$ is
countably incomplete, and therefore free. Thus, by
Theorem~(\ref{T: Ultrafilter existence}), there exists an
ultrafilter $\mathcal{U}$ on $\mathcal{D}(\mathbb{R}^{d})$
containing $\mathcal{F}_{\mathcal{B}}$. We shall keep
$\mathcal{U}$ fixed in what follows.

    \section{Non-Standard Numbers}\label{S:
Non-Standard Numbers} We will now use the ultrafilter defined in
the previous section to construct fields of non-standard real and
complex numbers.
\begin{definition}\label{D: Non-Standard Numbers}
Let $\mathcal{U}$ be as before, and let
$\mathbb{C}^{\mathcal{D}(\mathbb{R}^{d})}$ be the ring of
functionals from $\mathcal{D}(\mathbb{R}^{d})$ to $\mathbb{C}$
supplied with pointwise addition and multiplication. We shall
denote these functionals as ``families'' $(A_{\varphi})$ and treat
the domain $\mathcal{D}(\mathbb{R}^{d})$ as an ``index
set''.\\
We may define the operations of absolute value, real part mapping,
imaginary part mapping, and complex conjugation on the elements of
$\mathbb{C}^{\mathcal{D}(\mathbb{R}^{d})}$ by:
\begin{align*}
    &\vert (A_{\varphi})\vert = (\vert
A_{\varphi}\vert)\\
    &\Re (A_{\varphi}) = (\Re A_{\varphi})\\
    &\Im (A_{\varphi}) = (\Im A_{\varphi})\\
    &\overline{(A_{\varphi})} =
(\overline{A_{\varphi}})\\
\end{align*}
Also, we may define an embedding of  $\mathbb{C}$ into
$\mathbb{C}^{\mathcal{D}(\mathbb{R}^{d})}$ by $c \rightarrow
(C_{\varphi})$ where $C_{\varphi}=c$ for all $\varphi \in
\mathcal{D}(\mathbb{R}^{d})$.\\
Define an equivalence relation $\sim_{\mathcal{U}}$ on
$\mathbb{C}^{\mathcal{D}(\mathbb{R}^{d})}$ by
\[
    (A_{\varphi})\sim_{\mathcal{U}}(B_{\varphi}) \;\;
if \;\;
    \{\varphi\in\mathcal{D}(\mathbb{R}^{d}):
    A_{\varphi}=B_{\varphi}\}\in\mathcal{U}
    \]
Finally, let
$^{*}\mathbb{C}=\mathbb{C}^{\mathcal{D}(\mathbb{R}^{d})}/\sim_{\mathcal{U}}$.
That is, $^{*}\mathbb{C}$ consists of equivalence classes of
functionals in $\mathbb{C}^{\mathcal{D}(\mathbb{R}^{d})}$. We may
write $\langle (A_{\varphi})\rangle$ to represent these classes,
but to simplify notation we will denote by $\langle
A_{\varphi}\rangle\in$ $^{*}\mathbb{C}$ the non-standard number
(equivalence class of functionals) with representative
$(A_{\varphi})$.  $^{*}\mathbb{C}$ is called a field of complex
non-standard numbers.\\
$^{*}\mathbb{C}$ inherits the operations and embedding mentioned
above from $\mathbb{C}^{\mathcal{D}(\mathbb{R}^{d})}$. With the
embedding in mind, we shall treat elements of $\mathbb{C}$
as their images in $^{*}\mathbb{C}$.\\
A non-standard number $\langle A_{\varphi}\rangle$ is called
\textbf{real} if \\
\[
    \{\varphi\in \mathcal{D}(\mathbb{R}^{d}):
A_{\varphi} \in
    \mathbb{R}\}\in \mathcal{U}
    \]
We denote the set of all \textbf{real non-standard numbers} by
$^{*}\mathbb{R}$ and supply it with an order relation as follows:
\[
    \langle A_{\varphi}\rangle >_{ ^{*}\mathbb{R}} 0
\;\text{  if
    }\;
    \{\varphi \in \mathcal{D}(\mathbb{R}^{d}) :
    A_{\varphi}>0\}\in\mathcal{U}
    \]

\end{definition}

\begin{theorem}\label{T: Field Properties}
\begin{itemize}
    \item[(i)] Every number $\gamma\in\
^{*}\mathbb{C}$ can be
    uniquely represented in the form $\gamma = \alpha
+\beta i$
    where $\alpha, \beta\in$ $^{*}\mathbb{R}$ and
    $\alpha=\Re\gamma$, $\beta=\Im\gamma$, and
    $\vert\gamma\vert=\sqrt{\alpha^{2}+\beta^{2}}$.
    \item[(ii)] $^{*}\mathbb{C}$ is an algebraically
closed
    non-Archimedean field of characteristic zero.
$\mathbb{C}$ is
    a subfield of $^{*}\mathbb{C}$.
    \item[(iii)] $^{*}\mathbb{R}$ is a totally ordered
non-Archimedean
    real closed field. Moreover, $\alpha > 0$ in
$^{*}\mathbb{R}$
    \ifff $\alpha=\beta^{2}$ for some $\beta\in$
$^{*}\mathbb{R}$,
    $\beta\neq 0$. $\mathbb{R}$ is an ordered subfield
of
    $^{*}\mathbb{R}$.
\end{itemize}

\end{theorem}
\emph{Proof} (i) Let $\gamma\in\,^{*}\mathbb{C}$. Then
$\gamma=\langle C_{\varphi}\rangle$ for some
$(C_{\varphi})\in\mathbb{C}^{\mathcal{D}(\mathbb{R}^{d})}$. But
for each $\varphi$, $C_{\varphi}=A_{\varphi}+B_{\varphi}i$, where
$A_{\varphi}=\Re C_{\varphi}$ and $B_{\varphi}=\Im C_{\varphi}$.
Thus $\langle C_{\varphi}\rangle=\langle A_{\varphi}\rangle +
\langle B_{\varphi}\rangle i$. To prove uniqueness, suppose that
$\langle C_{\varphi}\rangle=\langle D_{\varphi}\rangle + \langle
E_{\varphi}\rangle i$ also. Then
\[
    \{\varphi :
    \Re C_{\varphi}=D_{\varphi}\}\cap\{\varphi :
    \Re C_{\varphi}=A_{\varphi}\}=\{\varphi :
A_{\varphi}=D_{\varphi}\}
    \in\mathcal{U}
    \]
    because $\mathcal{U}$ is closed under
intersections.
    Therefore $\langle A_{\varphi}\rangle=\langle
    D_{\varphi}\rangle$.\\
    The same argument can be applied to show that
$\langle B_{\varphi}\rangle=\langle
    E_{\varphi}\rangle$.\\
    The proof for $\vert\gamma\vert$ is similar.\\\\
(ii) It is not hard to check that
$\mathbb{C}^{\mathcal{D}(\mathbb{R}^{d})}$ really is a ring, and
that $\sim_{\mathcal{U}}$ really is an equivalence relation. It
follows that $^{*}\mathbb{C}$ is a (commutative) ring. To prove
that $^{*}\mathbb{C}$ is a field, we must show that each non-zero
element has a multiplicative inverse. For any non-zero
$\gamma\in\,^{*}\mathbb{C}$, we may choose a representative
$(C_{\varphi})$ such that $C_{\varphi}\ne 0$ for all $\varphi$.
Let $D_{\varphi}=1/C_{\varphi}$ and $\delta=\langle
D_{\varphi}\rangle$. Then $\delta\gamma=\langle 1
\rangle$.\\
Let
\[
    P(x)=\sum_{k=0}^{n}\alpha_{k}x^{k}, \;
\alpha_{k}\in\,
    ^{*}\mathbb{C} \text{ for all k }
    \]
    be a polynomial in $^{*}\mathbb{C}[x]$. Define
\[
    P_{\varphi}(x)= \sum_{k=0}^{n} A_{k,\varphi} x^{k}
    \]
    where $\alpha_{k}=\langle A_{k,\varphi}\rangle$
for each k. Since each $P_{\varphi}(x)$ is a polynomial over
$\mathbb{C}$, there exists a number $C_{\varphi}\in\mathbb{C}$
such that $P_{\varphi}(C_{\varphi})=0$. If we let $\gamma=\langle
C_{\varphi} \rangle$, it follows that $P(\gamma)=0$ in
$^{*}\mathbb{C}$.\\
That $\mathbb{C}$ is a subfield of $^{*}\mathbb{C}$ is clear from
the embedding.\\\\
(iii) The trichotomy of the order relation on $^{*}\mathbb{R}$
follows from the trichotomy of the order relation on $\mathbb{R}$.
For suppose $\mathcal{A}=\{\varphi : A_{\varphi}<B_{\varphi}\}$,
$\mathcal{B}=\{\varphi : A_{\varphi}=B_{\varphi}\}$, and
$\mathcal{C}=\{\varphi : A_{\varphi}>B_{\varphi}\}$, for some
non-standard real numbers $\langle A_{\varphi}\rangle$, $\langle
B_\varphi \rangle$. Note that $\mathcal{A}$, $\mathcal{B}$, and
$\mathcal{C}$ are mutually disjoint. Therefore, at most one of
$\mathcal{A}$, $\mathcal{B}$, or $\mathcal{C}$ can be in
$\mathcal{U}$. Also,
$\mathcal{A}\cup\mathcal{B}\cup\mathcal{C}=\mathcal{D}(\mathbb{R}^{d})\in\mathcal{U}$.
We can use this to prove that one of $\mathcal{A}$, $\mathcal{B}$,
or $\mathcal{C}$ must be in $\mathcal{U}$. For suppose that none
of $\mathcal{A}$, $\mathcal{B}$, or $\mathcal{C}$ is in
$\mathcal{U}$. Then by Theorem~(\ref{T: Ultrafilter
characterization}), $\mathcal{B}\cup\mathcal{C}\in\mathcal{U}$ and
$\mathcal{A}\cup\mathcal{C}\in\mathcal{U}$. Taking the
intersection of these two sets, we would have
$\mathcal{C}\in\mathcal{U}$, a contradiction. $\blacktriangle$

\begin{definition}\label{D: I(*C), F(*C), L(*C)}
Define the sets of \textbf{infinitesimal}, \textbf{finite}, and
\textbf{infinitely large} numbers as follows:
\[
    \mathcal{I}(^{*}\mathbb{C})=\{x\in\,
^{*}\mathbb{C} : \vert
    x\vert<1/n \text{ for all } n\in\mathbb{N}\}
    \]
\[
    \mathcal{F}(^{*}\mathbb{C})=\{x\in\,
^{*}\mathbb{C} : \vert
    x\vert<n \text{ for some } n\in\mathbb{N}\}
    \]
\[
    \mathcal{L}(^{*}\mathbb{C})=\{x\in\,
^{*}\mathbb{C} : \vert
    x\vert>n \text{ for all } n\in\mathbb{N}\}
    \]

\end{definition}
It is not hard to prove that $\mathcal{F}(^{*}\mathbb{C})$ is a
subring of $^{*}\mathbb{C}$ and $\mathcal{I}(^{*}\mathbb{C})$ is a
maximal ideal in $\mathcal{F}(^{*}\mathbb{C})$.

\begin{example}\label{E: Canonical infinitesimal}
Define $(R_{\varphi})\in\mathbb{C}^{\mathcal{D}(\mathbb{R}^{d})}$
by
\[
    R_{\varphi}= \sup\{\Vert x\Vert : x\in \supp
\varphi\}
    \]
    where $\supp
\varphi=\overline{\{x\in\mathbb{R}^{d} : \varphi(x)\ne
    0\}}$ is the support of $\varphi$.
    The non-standard number
    $\rho=\langle R_{\varphi}\rangle$ is a (positive)
infinitesimal.
     For let $\mathcal{A}=\{\varphi :
    0< R_{\varphi} < 1/n\}$. Then for any
    $\varphi\in\mathcal{B}_{n+1}$, we have
    $\varphi\in\mathcal{A}$, by the definition of
    $\mathcal{B}_{n+1}$. Thus
    $\mathcal{B}_{n+1}\subset\mathcal{A}$, implying
    $\mathcal{A}\in\mathcal{U}$. $\rho$ is called the
canonical
    infinitesimal in  $^{*}\mathbb{C}$.
\end{example}

\begin{definition}\label{D: standard part mapping}
Define the \textbf{standard part mapping} $st:$
$^{*}\mathbb{R}\rightarrow\mathbb{R}\cup\{\pm\infty\}$ by
\[
    st(x)=
    \begin{cases}
        \sup\{r\in\mathbb{R} : r<x\} & \text{if
        $x\in\mathcal{F}(^{*}\mathbb{R})$}\\
        \infty & \text{if
$x\in\mathcal{L}(^{*}\mathbb{R}_{+})$}\\
        -\infty & \text{if
$x\in\mathcal{L}(^{*}\mathbb{R}_{-})$}
    \end{cases}
    \]
We may extend this definition to $^{*}\mathbb{C}$ by
$st(x+yi)=st(x)+st(y)i$.
\end{definition}

\begin{theorem}\label{T: asymptotic expansion}
If $x\in\mathcal{F}(^{*}\mathbb{C})$ then $x$ has a unique
asymptotic expansion: $x=r + dx$ where $r\in\mathbb{C}$ and
$dx\in\mathcal{I}(^{*}\mathbb{C})$. In fact, $r=st(x)$.
\end{theorem}
\emph{Proof} We will prove the case for
$x\in\mathcal{F}(^{*}\mathbb{R})$. The general result will follow.
Let $x\in\mathcal{F}(^{*}\mathbb{R})$. First note that
$x-st(x)\in\mathcal{I}(^{*}\mathbb{R})$, for otherwise we would
have $\vert x-st(x)\vert>1/n$ for some n, implying either that
$st(x)>x$ or that $st(x)+1/2n<x$. In either case, this is a
contradiction to Definition~(\ref{D: standard part mapping}). To
prove uniqueness, suppose that $x=r+dx$ and $x=s+dy$ are two
expansions of $x$. Then we would have $r-s=dx-dy$, implying that
$r-s\in\mathcal{I}(^{*}\mathbb{R})$. But since $r-s\in\mathbb{R}$,
$r-s=0$. Hence $r=s$. Therefore $r+dx=r+dy$, implying $dx=dy$.
$\blacktriangle$\\

    \section{Internal Sets}\label{S: Internal Sets}
In non-standard analysis, internal sets play the role of the
``good" sets, in a similar way to the measurable sets in Lebesgue
theory.\\
In what follows we will use the abbreviation \textbf{a.e.} to mean
that the set of functions for which some statement is true is in
$\mathcal{U}$.

\begin{definition}\label{D: non-standard extension}
Let $\mathbb{A}\subseteq\mathbb{C}$. The \textbf{non-standard
extension} of $\mathbb{A}$ is
\[
    ^{*}\mathbb{A}=\{\langle A_{\varphi}\rangle\in\,
^{*}\mathbb{C} :
    A_{\varphi}\in\mathbb{A}\; a.e.\}
    \]
A set $\mathcal{A}$ of non-standard numbers is called
\textbf{internal standard} if it is the non-standard extension of
some subset of $\mathbb{C}$. The set of all internal standard sets
is denoted by  $^{\sigma}\mathcal{P}(\mathbb{C})$.

\end{definition}

\begin{example}\label{E: internal standard intervals}
The non-standard extensions of the intervals $(a,\,b)$, $[a,b]$,
$(a, \infty)$, etc. are
\begin{align*}
    &^{*}(a, b) = \{x\in\ ^{*}\mathbb{R} : \,a<x<b\}\\
    &^{*}[a, b] = \{x\in\ ^{*}\mathbb{R} : \,a\leq
x\leq b\}\\
    &^{*}(a, \infty) = \{x\in\ ^{*}\mathbb{R} :
\,a<x\},\;etc.
\end{align*}
\end{example}

\begin{definition}\label{D: internal sets}
Let
$(\mathbb{A}_{\varphi})\in\mathcal{P}(\mathbb{C})^{\mathcal{D}(\mathbb{R}^{d})}$
be a family of subsets of $\mathbb{C}$. We define the
\textbf{internal set} generated by $(\mathbb{A}_{\varphi})$ by
\[
    \langle\mathbb{A}_{\varphi}\rangle=\{\langle
    A_{\varphi}\rangle\in\, ^{*}\mathbb{C} :
    A_{\varphi}\in\mathbb{A}_{\varphi} \; a.e.\}
\]
A set is called \textbf{external} if it is not internal.
\end{definition}

\begin{example}\label{E: internal set (0,rho)}
Let $\mathbb{A}_{\varphi}=(0,R_{\varphi})$, where $R_{\varphi}$ is
as in Example~(\ref{E: Canonical infinitesimal}). Then the
internal set $\langle \mathbb{A}_{\varphi} \rangle$ generated by
$(\mathbb{A}_{\varphi})$ is the internal interval $(0, \rho)$. It
is important to note (and easy to check) that this coincides with
the more natural definition for $(0, \rho)$ given by
\[
    (0, \rho) = \{x\in\, ^{*}\mathbb{R} : 0 < x <
\rho\}
    \]
\end{example}
A set $S\subset\mathbb{R}^d$ is called {\em relatively compact} if
its closure $\overline{S}$ is compact in $\mathbb{R}^d$. Unless it
is specified otherwise, we shall call {\em Lebesgue measurable
sets} of $\mathbb{R}^d$ simply {\em measurable} sets.

\begin{definition}\label{D: *-Properties} An internal
set $\bra\mathbb{A}_\varphi\ket$ of $^*\mathbb{R}^d$ is called
$*$-{\bf measurable ($*$-compact, $*$-relatively-compact,
$*$-closed, $*$-open, etc.)} if
\[
    \mathbb{A}_{\varphi} \text{\;is measurable
(compact, relatively compact, closed, etc.) in \;} \mathbb{R}^{d}
\text{ for a.e.
 } \varphi
    \]

\end{definition}
Let $\rho\in{^*\mathbb{R}}$ denote a positive infinitesimal in
$^*\mathbb{R}$ (for example, $\rho$ might be the positive
infinitesimal defined in (Example~\ref{E: Canonical
infinitesimal})). We shall keep $\rho$ fixed in what follows.

\begin{definition}\label{D: rho-Moderate ad rho-Finite
Numbers} Let $\rho$ be a positive infinitesimal in $^*\mathbb{R}$.
We define the following (external) sets of non-standard numbers:
\begin{align}
&\mathcal{M}_\rho({^*\mathbb{C}})=\{x\in{^*\mathbb{C}}\mid \;
|x|\leq\rho^{-n} \text{\: for some\;}
n\in\mathbb{N}\}\notag\\
&\mathcal{N}_\rho({^*\mathbb{C}})=\{x\in{^*\mathbb{C}}\mid \;
|x|<\rho^{n} \text{\: for all\;}
n\in\mathbb{N}\}\notag \\
&\mathcal{F}_\rho(^*\mathbb{C})=\{x\in{^*\mathbb{C}}\mid \;
|x|<1/\sqrt[n]{\rho} \text{\; for all\; }
n\in\mathbb{N}\},\notag\\
&\mathcal{I}_\rho(^*\mathbb{C})=\{x\in{^*\mathbb{C}}\mid \;
|x|\leq\sqrt[n]{\rho} \text{\; for some\; }
n\in\mathbb{N}\},\notag\\
&\mathcal{C}_\rho(^*\mathbb{C})=\{x\in{^*\mathbb{C}}\mid \;
\sqrt[n]{\rho}<|x|<1/\sqrt[n]{\rho} \text{\; for all\; }
n\in\mathbb{N}\}.\notag
\end{align}
The numbers in $\mathcal{M}_\rho(^*\mathbb{C})$ and
$\mathcal{N}_\rho(^*\mathbb{C})$ are called
\textbf{$\rho$-moderate} and \textbf{$\rho$-null} non-standard
numbers, respectively. Similarly, the numbers in
$\mathcal{F}_\rho(^*\mathbb{C})$, $\mathcal{I}_\rho(^*\mathbb{C})$
and $\mathcal{C}_\rho(^*\mathbb{C})$ are called
\textbf{$\rho$-finite, $\rho$-infinitesimal} and
\textbf{$\rho$-constant}, respectively.

\end{definition}

    \section{Saturation Principle in
$^*\mathbb{C}$}\label{S: Saturation Principle in *C}
\begin{theorem}\label{T: saturation principle}
Let $\{\mathcal{A}_{n}\}$ be a sequence of internal sets in
$^{*}\mathbb{C}$ such that
\[
    \bigcap_{n=0}^{m}\mathcal{A}_{n}\ne\varnothing
    \]
    for all $m\in\mathbb{N}$. (The sequence
$\{\mathcal{A}_{n}\}$
    satisfies the finite intersection property.) Then
\[
\bigcap_{n=0}^{\infty}\mathcal{A}_{n}\ne\varnothing.
    \]
\end{theorem}
\emph{Proof} Since each $\mathcal{A}_{n}$ is internal,
\[
    \mathcal{A}_{n}=\langle
\mathbb{A}_{n,\varphi}\rangle \text{ ,
    $\mathbb{A}_{n, \varphi}\subseteq\mathbb{C}$}
    \]
Also, since for each $m$ it is given that
$\bigcap_{n=0}^{m}\mathcal{A}_{n}\ne\varnothing$, this implies
that for each $m$ there exists a non-standard number $\langle
C_{m, \varphi}\rangle\in$ $^{*}\mathbb{C}$ such that
\[
    \langle C_{m,\varphi} \rangle \in
\bigcap_{n=0}^{m}
    \langle \mathbb{A}_{n,\varphi}\rangle
    \]
or, in other words,
\[
    \langle C_{m,\varphi} \rangle \in \langle
    \mathbb{A}_{n,\varphi} \rangle \text{ for } 0\leq
n \leq m
    \]
This means that for a.e. $\varphi$ and $0 \leq n \leq m$,
\[
    C_{m,\varphi} \in \mathbb{A}_{n,\varphi}
    \]
Remembering that $\mathcal{U}$ is closed under finite
intersections, we see that for a.e. $\varphi$,
\[
    C_{m,\varphi} \in
\bigcap_{n=0}^{m}\mathbb{A}_{n,\varphi}
    \]
Hence for each $m$,
\[
    \bigcap_{n=0}^{m}\mathbb{A}_{n,\varphi} \ne
\varnothing \; a.e.
    \]
We may assume without loss of generality that
$\mathbb{A}_{0,\varphi}$ is non-empty for all $\varphi$. (Else
define $\mathbb{A'}_{0,\varphi}=\mathbb{A}_{0,\varphi}$ if
$\mathbb{A}_{0,\varphi}\ne\varnothing$ and
$\mathbb{A'}_{0,\varphi}=\mathbb{C}$ otherwise. Then it will still
be true that $\mathcal{A}_{0}=\langle \mathbb{A'}_{0,\varphi}
\rangle$ .)\\
Next define a function $\mu : \mathcal{D}(\mathbb{R}^{D})
\longrightarrow\mathbb{N}\cup \{ \infty\}$ by
\[
\mu(\varphi)=\max\{m\in\mathbb{N}_{0}\cup\{\infty\} \, \vert\,
\bigcap_{n=0}^{m}\mathbb{A}_{n,\varphi}\ne\varnothing\}
    \]
Notice that $\mu$ is defined for all $\varphi$ due to our
assumption for $\mathbb{A}_{0,\varphi}$.\\
Thus we have
\[
\bigcap_{n=0}^{\mu(\varphi)}\mathbb{A}_{n,\varphi}\ne\varnothing
    \text{ for all }
\varphi\in\mathcal{D}(\mathbb{R}^{d})
    \]
Hence for every $\varphi\in\mathcal{D}(\mathbb{R}^{d})$ there
exists (by Axiom of Choice) $A_{\varphi}$ such that
$A_{\varphi}\in\bigcap_{n=0}^{\mu(\varphi)}\mathbb{A}_{n,\varphi}$.\\
We intend to show that
\[
    \langle A_{\varphi} \rangle \in
    \bigcap_{n=0}^{\infty}\mathcal{A}_{n}
    \]
or, equivalently, that for every $m$,
$A_{\varphi}\in\mathbb{A}_{m,\varphi}$ for a.e.
$\varphi$.\\
If $\varphi$ is such that
$\bigcap_{n=0}^{m}\mathbb{A}_{n,\varphi}\ne\varnothing$, this
implies that $0\leq m\leq\mu(\varphi)$. Thus
$A_{\varphi}\in\mathbb{A}_{m,\varphi}$, by the choice of
$A_\varphi$. Therefore,
\[
    \{\varphi \, \vert \,
\bigcap_{n=0}^{m}\mathbb{A}_{n,\varphi}\ne\varnothing\}
     \subseteq \{\varphi \, \vert \,
    A_{\varphi}\in\mathbb{A}_{m,\varphi}\}
    \]
But the set on the left is in $\mathcal{U}$, and so the set on the
right is also, as required. $\blacktriangle$ \\

    \section{Non-Standard Smooth Functions}\label{S:
Non-Standard Smooth Functions}

Having constructed the fields $^{*}\mathbb{R}$ and $^*\mathbb{C}$,
the natural next step is to look at functions on these fields.
However, for our purposes we will focus on a certain class of
function contained in $^{*}\mathbb{C}^{^{*}\mathbb{R}}$. In what
follows $\mathcal{E}(\mathbb{R})$ is the set of
$\mathcal{C}^{\infty}$-functions from $\mathbb{R}$ into
$\mathbb{C}$.

\begin{definition}\label{D: non-standard smooth
functions} A function $f\in\,^{*}\mathbb{C}^{^{*}\mathbb{R}}$ is
called \textbf{internal smooth} if there exists a family
$(f_{\varphi})\in\mathcal{E}(\mathbb{R})^{\mathcal{D}(\mathbb{R}^{d})}$
such that for every $x=\langle X_{\varphi} \rangle \in \,
^{*}\mathbb{R}$
\[
    f(x) = \langle f_{\varphi} (X_\varphi) \rangle
    \]
The set of all internal smooth functions will be denoted by
$^{*}\mathcal{E}(\mathbb{R})$.
\end{definition}

\begin{remark}\label{R: equivalent def of *E }
$^{*}\mathcal{E}(\mathbb{R})$ may equivalently be defined as the
set of equivalence classes $\langle f_{\varphi} \rangle $ of
families of functions in
$\mathcal{E}(\mathbb{R})^{\mathcal{D}(\mathbb{R}^{d})}$, where the
equivalence relation is as usual:
\[
    (f_{\varphi}) \sim_{\mathcal{U}} (g_{\varphi})
\text{ if }
    f_{\varphi} = g_{\varphi} \text{ for a.e. }
\varphi
    \]
    It is not hard to prove that the value
    of an internal function does not depend on the
    choice of representatives. If $\langle X_{\varphi}
\rangle
    = \langle Y_{\varphi} \rangle \in \,
^{*}\mathbb{R}$ and $\langle f_{\varphi} \rangle
    = \langle g_{\varphi} \rangle \in \,
    ^{*}\mathcal{E}(\mathbb{R})$, then
    \[
        \{\varphi \, \vert \,
f_{\varphi}(X_\varphi)=g_{\varphi}(X_\varphi)\} \cap \{\varphi \,
\vert \,
g_{\varphi}(X_\varphi)=g_{\varphi}(Y_\varphi)\}\subseteq\{\varphi
\, \vert \, f_{\varphi}(X_\varphi)=g_{\varphi}(Y_\varphi)\}
        \]
    Since $\mathcal{U}$ is closed under intersections,
$\langle
    f_{\varphi} (X_{\varphi}) \rangle = \langle
g_{\varphi}
    (Y_{\varphi}) \rangle$.
\end{remark}
The operations of addition, multiplication, and partial
differentiation in $^{*}\mathcal{E}(\mathbb{R}^{d})$ are inherited
from $\mathcal{E}(\mathbb{R}^{d})$. Also,
$\mathcal{E}(\mathbb{R}^{d})$ is embedded in
$^{*}\mathcal{E}(\mathbb{R}^{d})$ by $f \longrightarrow \, ^{*}f$
where $^{*}f=\langle f_{\varphi} \rangle$, $f_{\varphi}=f$ for all
$\varphi \in \mathcal{D}(\mathbb{R}^{d})$.\\
In what follows {\em integrable} means {\em Lebesgue integrable}.

\begin{definition}\label{D: Integration in *E}  Let
$\bra\mathbb{X}_\varphi\ket\subseteq \, ^*\mathbb{R}^d$ be a
$*$-measurable internal set and let $\bra
f_\varphi\ket\in{^*\mathcal{E}(\mathbb{R}^d)}$ be an internal
function. We say that $\bra f_\varphi\ket$ is {\bf $*$-integrable}
over $\bra\mathbb{X}_\varphi\ket$ if
\[
     f_\varphi \text{\; is integrable over\; }
\mathbb{X}_\varphi \text{ for a.e. } \varphi
    \]
If $\bra f_\varphi\ket$ is {\bf $*$-integrable} over
$\bra\mathbb{X}_\varphi\ket$, we define the integral:
\begin{equation}\label{E: Integral}
\int_{\bra\mathbb{X}_\varphi\ket}\bra f_\varphi\ket(x)\, dx=\left<
\int_{\mathbb{X}_\varphi} f_\varphi(x)\, dx\right>.
\end{equation}
We also say that the integral {\bf converges} in $^*\mathbb{C}$
(since it is a number in $^*\mathbb{C}$). Notice that as long as
the integral converges for a.e. $\varphi$, we include this object
in the equivalence class, even if the integral diverges for other
$\varphi$.
\end{definition}

    \section{Internal Sets and Saturation Principle in
$^* \mathcal{E}(\mathbb{R}^{d})$}
    \label{S: Internal Sets and Saturation Principle
in *E} We define internal sets in $^*\mathcal{E}(\mathbb{R}^{d})$
similarly to those of $^*\mathbb{C}$.

\begin{definition}\label{D: internal sets in *E}
\begin{itemize}
    \item[(i)] Let
$(\mathbb{F}_{\varphi})\in\mathcal{P}(\mathcal{E}(\mathbb{R}^{d}))
^{\mathcal{D}(\mathbb{R}^{d})}$ be a family of subsets of
$\mathcal{E}(\mathbb{R}^{d})$. We define the \textbf{internal set}
generated by $(\mathbb{F}_{\varphi})$ by
\[
    \langle\mathbb{F}_{\varphi}\rangle=\{\langle
    f_{\varphi}\rangle\in\,
^{*}\mathcal{E}(\mathbb{R}^{d}) :
    f_{\varphi}\in\mathbb{F}_{\varphi} \; a.e.\}
\]
A set is called \textbf{external} if it is not internal.

    \item[(ii)] An internal set $\mathcal{F}$ is
called
    \textbf{standard} if there exists
$\mathbb{F}\subset
    \mathcal{E}(\mathbb{R}^{d})$ such that
$\mathcal{F}=\langle
    \mathbb{F} \rangle$. In this case we may also
write
    $\mathcal{F}= \, ^*\mathbb{F}$.

\end{itemize}

\end{definition}

\begin{theorem}\label{T: saturation principle in *E}
Let $\{\mathcal{F}_{n}\}$ be a sequence of internal sets in
$^{*}\mathcal{E}(\mathbb{R}^{d})$ such that
\[
    \bigcap_{n=0}^{m}\mathcal{F}_{n}\ne\varnothing
    \]
    for all $m\in\mathbb{N}$. (The sequence
$\{\mathcal{F}_{n}\}$
    satisfies the finite intersection property.) Then
\[
\bigcap_{n=0}^{\infty}\mathcal{F}_{n}\ne\varnothing.
    \]
\end{theorem}
\emph{Proof} The proof is almost identical to that of
(Theorem~\ref{T: saturation principle}).

\begin{definition}\label{D: rho-Moderate and
rho-Finite Functions} We define the following (external) subsets
of ${^*\mathcal{E}(\mathbb{R}^{d})}$:
\begin{align}\notag
&\mathcal{F}(^*\mathcal{E}(\mathbb{R}^{d}))=\{f\in{^*\mathcal{E}(\mathbb{R}^{d})}\mid
(\forall\alpha\in\mathbb{N}_0^d)(\forall
x\in\mathcal{F}(^*\mathbb{R}^{d}))\left[\partial^\alpha
f(x)\in\mathcal{F}(^*\mathbb{C})\right]\},\\\notag
&\mathcal{I}(^*\mathcal{E}(\mathbb{R}^{d}))=\{f\in{^*\mathcal{E}(\mathbb{R}^{d})}\mid
(\forall\alpha\in\mathbb{N}_0^d)(\forall
x\in\mathcal{F}(^*\mathbb{R}^{d}))\left[\partial^\alpha
f(x)\in\mathcal{I}(^*\mathbb{C})\right]\},\notag\\
&\mathcal{M}_\rho(^*\mathcal{E}(\mathbb{R}^{d}))=\left\{f\in{^*\mathcal{E}(\mathbb{R}^{d})}
\mid(\forall\alpha\in\mathbb{N}_0^d)(\forall
x\in\mathcal{F}(^*\mathbb{R}^{d}) )\left[\partial^\alpha
f(x)\in\mathcal{M}_\rho(^*\mathbb{C})\right]\right\},\notag\\
&\mathcal{N}_\rho(^*\mathcal{E}(\mathbb{R}^{d}))=\left\{f\in{^*\mathcal{E}(\mathbb{R}^{d})}
\mid(\forall\alpha\in\mathbb{N}_0^d)(\forall
x\in\mathcal{F}(^*\mathbb{R}^{d}))\left[\partial^\alpha
f(x)\in\mathcal{N}_\rho(^*\mathbb{C})\right]\right\},\notag\\
&\mathcal{F}_\rho(^*\mathcal{E}(\mathbb{R}^{d}))=\left\{f\in{^*\mathcal{E}(\mathbb{R}^{d})}
\mid(\forall\alpha\in\mathbb{N}_0^d)(\forall
x\in\mathcal{F}(^*\mathbb{R}^{d}))\left[\partial^\alpha
f(x)\in\mathcal{F}_\rho(^*\mathbb{C})\right]\right\},\notag\\
&\mathcal{I}_\rho(^*\mathcal{E}(\mathbb{R}^{d}))=\left\{f\in{^*\mathcal{E}(\mathbb{R}^{d})}
\mid(\forall\alpha\in\mathbb{N}_0^d)(\forall
x\in\mathcal{F}(^*\mathbb{R}^{d}))\left[\partial^\alpha
f(x)\in\mathcal{I}_\rho(^*\mathbb{C})\right]\right\},\notag\\
&\mathcal{C}_\rho(^*\mathcal{E}(\mathbb{R}^{d}))=\left\{f\in{^*\mathcal{E}(\mathbb{R}^{d})}
\mid(\forall\alpha\in\mathbb{N}_0^d)(\forall
x\in\mathcal{F}(^*\mathbb{R}^{d}))\left[\partial^\alpha
f(x)\in\mathcal{C}_\rho(^*\mathbb{C})\right]\right\}.\notag
\end{align}
The functions in $\mathcal{F}(^*\mathcal{E}(\mathbb{R}^{d})),
\mathcal{I}(^*\mathcal{E}(\mathbb{R}^{d}))$,
$\mathcal{M}_\rho(^*\mathcal{E}(\mathbb{R}^{d}))$,
$\mathcal{N}_\rho(^*\mathcal{E}(\mathbb{R}^{d})$,
$\mathcal{F}_\rho(^*\mathcal{E}(\mathbb{R}^{d})$,
$\mathcal{I}_\rho(^*\mathcal{E}(\mathbb{R}^{d})$, and
$\mathcal{C}_\rho(^*\mathcal{E}(\mathbb{R}^{d})$ are called
\textbf{ finite, infinitesimal, $\rho$-moderate, $\rho$-null,
$\rho$- finite, $\rho$-infinitesimal and $\rho$-constant
functions}, respectively. For more details we refer to (Lightstone
and Robinson~\cite{LiRob}) and (Wolf and Todorov~\cite{TodWolf}.

\end{definition}

    \section{Weak Equality}\label{S: Weak equality}

\begin{definition}\label{D: weak equalities}
Let $x,y\in\, ^{*}\mathbb{C}$, $f,\, g \in \,
^{*}\mathcal{E}(\mathbb{R}^{d})$
    \begin{itemize}
        \item[(i)] $x \approx y$ if $ x-y \in
                    \mathcal{I}(^{*}\mathbb{C})$
        \item[(ii)] $x \overset{\rho}{=} y$ if $ x - y
                    \in
\mathcal{N}_{\rho}(^{*}\mathbb{C})$
        \item[(iii)] $f \approx g$ if $f-g \in
\mathcal{I}(^{*}\mathcal{E}(\mathbb{R}^{d}))$
        \item[(iv)] $f \overset{\rho}{=} g$ if $f-g\in
\mathcal{N}_{\rho}(^{*}\mathcal{E}(\mathbb{R}^{d}))$

        \item[(v)] $f \cong g $ if $ \int
f(x)\tau(x)\, dx = \int
                    g(x) \tau (x) \, dx$ for every

$\tau\in\mathcal{D}(\mathbb{R}^{d})$
        \item[(vi)] $f \overset{\rho}{\simeq} g$ if $
\int f(x)\tau(x)\, dx
                    \overset{\rho}{=} \int
                    g(x) \tau (x) \, dx$ for every

$\tau\in\mathcal{D}(\mathbb{R}^{d})$

        \item[(vii)] $f \approxeq g$ if $ \int
f(x)\tau(x)\, dx \approx \int
                    g(x) \tau (x) \, dx$ for every

$\tau\in\mathcal{D}(\mathbb{R}^{d})$

    \end{itemize}

\end{definition}
It is not hard to prove that each of these weak equalities forms
an equivalence relation in its respective space. Many results in
non-standard analysis hold \textbf{weakly} in the sense of one of
these weak equalities.

    \section{Schwartz Distributions}\label{S: Schwartz
    Distributions}
    At this point, we must take a short detour to
present
    some basic definitions and results from the
Schwartz theory.

\begin{definition}\label{D: distribution}
A \textbf{distribution} is a mapping $F:
\mathcal{D}(\mathbb{R}^{d}) \longrightarrow \mathbb{C}$ that
satisfies the following conditions:
    \begin{itemize}
            \item[(i)] Linearity:

$F[c_{1}\tau_{1}+c_{2}\tau_{2}]=c_{1}F[\tau_{1}]+c_{2}F[\tau_{2}]$
            for all $c_{1},c_{2}\in\mathbb{C}$ and

$\tau_{1},\tau_{2}\in\mathcal{D}(\mathbb{R}^{d})$.
            \item[(ii)] Continuity: Let $\{\tau_{k}\}$
be a
            sequence in $\mathcal{D}(\mathbb{R}^{d})$.
Suppose there exists
            $R$ such that $\supp\tau_{k}\subseteq \{ x
: \; |x|<R
            \}$ for all $k$. Also, suppose there
exists
            $\tau\in\mathcal{D}(\mathbb{R}^{d})$ such
that
$\partial^{\alpha}\tau_{k}\longrightarrow\partial^{\alpha}\tau$
            for all $\alpha\in\mathbb{N}_{0}^{d}$
uniformly as
            $k\longrightarrow\infty$. Then
            $F[\tau_{k}]\longrightarrow F[\tau]$.

\end{itemize}
We will denote by $\mathcal{D}'(\mathbb{R}^{d})$ the set of all
such distributions.

\end{definition}
We supply $\mathcal{D}'(\mathbb{R}^{d})$ with the usual pointwise
addition and scalar multiplication. In addition, we define partial
differentiation by
\[
    (\partial^{\alpha}F)[\tau] = (-1)^{|\alpha|
    }F[\partial^{\alpha}\tau]
    \]
and multiplication by a smooth function $g \in
\mathcal{E}(\mathbb{R}^{d})$ by
\[
    (gF)[\tau]=F[g\tau]
    \]
    Both of these operations are well-defined since
    $\partial^{\alpha}\tau, \,
    g\tau\in\mathcal{D}(\mathbb{R}^{d})$.\\
$\mathcal{L}_{loc}(\mathbb{R}^{d})$, the set of locally integrable
functions, is embedded in $\mathcal{D}'(\mathbb{R}^{d})$ by the
mapping
\[
    S(f)=\int f(t) \tau (t) dt
    \]
It is not hard to show that this embedding preserves the
operations mentioned above.\\
Finally, we define the \textbf{convolution} of a distribution with
a test function by
\[
    (F \ast \tau)(x) = F[\tau(x-t)]
    \]

\begin{theorem}\label{T: F*phi in E}
If $F\in\mathcal{D}'(\mathbb{R}^{d})$ and
$\tau\in\mathcal{D}(\mathbb{R}^{d})$ then
$(F\ast\tau)(x)\in\mathcal{E}(\mathbb{R}^{d})$ and
$\partial^{\alpha}(F\ast\tau)=F\ast\partial^{\alpha}\tau$.

\end{theorem}
Before proving this theorem, we will state (without proof) a
result from analysis. See (Rudin~\cite{Rudin} p.148):

\begin{lemma}\label{L: uniform convergence}
Suppose
\[
    \lim_{n\to\infty} f_n(x) = f(x) \quad (x\in E).
    \]
Put
\[
    M_n = \sup_{x\in E} \vert f_n(x)-f(x) \vert .
    \]
 Then $f_n \longrightarrow f$ uniformly on $E$ if and
only if $M_n
 \longrightarrow 0$ as $n \longrightarrow \infty$ .

\end{lemma}
\emph{Proof of the theorem} We will prove the theorem for the case
d=1. The general result will follow.\\
Let $f(x) = (F\ast\tau)(x)$. Fixing $x$, we wish to show that
\[
    f(x+h)-f(x) \longrightarrow 0 \text{ as }
h\longrightarrow 0
    \]
Note that
\begin{align}\notag
    f(x+h)-f(x) &=
(F\ast\tau)(x+h)-(F\ast\tau)(x)\\\notag
                &= F[\tau(x+h-t)]-F[\tau(x-t)]\\\notag
                &= F[\tau(x+h-t)-\tau(x-t)]\notag
\end{align}
    by the linearity of $F$.
Let
\[
    \psi(t) = \tau(x+h-t)-\tau(x-t)
    \]
$\psi$ is itself a test function, and if we restrict $|h|<1$, then
the support of $\psi$ and all its derivatives is contained in $E =
\{y \, \vert \, |y| \leq r+ |x| +1\}$, where $r$ is the radius of
the support of $\tau$. It is clear that any sequence
$\{\psi_{h_{n}}\}$ where $h_n \longrightarrow 0$ as
$n\longrightarrow \infty$ converges pointwise to 0 for all $x$ and
$t$ (by the uniform continuity of $\tau$). Also, since one compact
set, $E$, contains the support of $\psi_{h_n}$ for all n, and
since each $\psi_{h_n}$ is continuous, $M_n = \sup_{t\in E}
|\psi_{h_n}(t)|$ is achieved by $\psi_{h_n}$ for each $n$. Thus
$M_n \to 0$, implying that $\{\psi_{h_{n}}\} \to 0$ uniformly, by
the Lemma. Therefore, since $F$ is continuous in the sense of
(Definition~\ref{D: distribution}),
\[
    F[\psi_{h_{n}}]\longrightarrow F[0] = 0
    \]
Since $x$ was chosen arbitrarily, this proves that $f=F\ast\tau$
is continuous.\\
To prove that $f'$ exists and that $f'(x) = (F\ast\tau ')(x)$, we
must show that
\[
    \frac{f(x+h) - f(x)}{h} - (F\ast\tau ')(x)
\longrightarrow 0
    \text{ as } h\longrightarrow 0
    \]
Note that
\begin{align}\notag
    \frac{f(x+h) - f(x)}{h} - (F\ast\tau ')(x)
    &= \frac{(F\ast\tau)(x+h)-(F\ast\tau)(x)}{h} -
(F\ast\tau')(x)\\\notag
    &= F\left[ \frac{\tau(x+h-t) - \tau(x-t)}{h}  -
\tau '(x-t)\right]\notag
\end{align}
Now, if we let
\[
    \chi(t) = \frac{\tau(x+h-t) - \tau(x-t)}{h}  -
\tau '(x-t)
    \]
we can use the same argument as before to show that
\[
    F[\chi_{h_{n}}(t)]\longrightarrow F[0]=0 \text{ as
}
    h\longrightarrow 0
    \]
This proves that $(F\ast\tau)' = F\ast\tau '$. Since $\tau '$ is
itself a test function, the same proof works to show that $(F\ast
\tau)'' = F\ast\tau ''$ and so on. For functions of several
variables, the same argument can be applied in each variable to
show the general result. $\blacktriangle$\\\\
Before we can prove the embedding of the distributions in
$^*\mathcal{E}(\mathbb{R}^{d})$, we need a result showing that
distributions can be ``approximated" in a way by a certain
sequence of test functions.\\
\begin{theorem}\label{T: regularization}
Let $\{\delta_n\}$ be a sequence in $\mathcal{D}(\mathbb{R}^{d})$
such that $\delta_n\in\mathcal{B}_n$ for every $n$. Then for any
distribution $T\in\mathcal{D}'(\mathbb{R}^{d})$, $T\ast\delta_n\to
T$ weakly. (A sequence of distributions $\{F_k\}$
\textbf{converges weakly} to a distribution $F$ if $F_k[\tau] \to
F[\tau]$ for all test functions $\tau$.)
\end{theorem}
Before proving this theorem, we need two lemmas:\\
\begin{lemma}\label{L: deltan*tau->tau}
Let $\{\delta_n\}$ be as above and let $\tau$ be any test
function. Then there exists $R$ such that
$\supp(\delta_n\ast\tau)\subset \{x : |x|\leq R\}$. Also,
$\partial^\alpha(\delta_n\ast\tau)\to\tau$ uniformly for every
$\alpha\in\mathbb{N}_0^d$.
\end{lemma}
\emph{Proof} For each $n$, $\supp\delta_n\subset\{x : |x|\leq
1/n\}$. In particular, $\supp\delta_n\subset\{x : |x| \leq 1\}$.
If we let $R_\tau$ be the radius of the support of $\tau$ and set
$R=R_\tau+1$, then it is not hard to see that
$\supp(\delta_n\ast\tau)\subset\{x : |x|\leq R\}$.\\
As before, to show the uniform convergence it is enough to prove
that
\[
    \sup_{|x|\leq R}|(\delta_n\ast\tau)(x)-\tau(x)|
\to 0
    \]
Recalling that $\int \delta_n =1$, we see that
\begin{align}\notag
    \sup_{|x|\leq R}|(\delta_n\ast\tau)(x)-\tau(x)| &=
    \sup_{|x|\leq R} \left| \int
\delta_n(t)\tau(x-t)dt -
    \tau(x)\int\delta_n(t)dt\right|\\\notag
    &=\sup_{|x|\leq R} \left| \int_{|t|\leq
    1/n}\delta_n(t)[\tau(x-t)-\tau(x)]dt\right|\notag
    \end{align}
By the mean value theorem for integrals, there exists $|t_n|\leq
1/n$ such that
\[
    =\sup_{|x|\leq
R}\left|[\tau(x-t_n)-\tau(x)]\int\delta_n(t)dt\right|
    \]
and by the extreme value theorem there exists $|x_n|\leq R$ such
that
\[
    =\left|\tau(x_n-t_n)-\tau(x_n)\right|
    \]
This last expression vanishes as $n\to\infty$ since $\tau$ is
uniformly continuous. The case $\alpha\neq 0$ is
similar. $\blacktriangle$\\\\
For the proof of the next lemma see (Folland~\cite{Folland}
p.318):
\begin{lemma}\label{L: Folland lemma}
Suppose $F$ is a distribution and $\phi$ and $\psi$ are test
functions. Then $(F\ast\phi)[\psi]=F[\tilde{\phi}\ast\psi]$, where
$\tilde{\phi}(x)=\phi(-x)$.
\end{lemma}
\emph{Proof of the theorem} We must show that for any distribution
$T$ and any test function $\tau$,
\[
    (T\ast\delta_n)[\tau(x)]\to T[\tau(x)]
    \]
Using (Lemma~\ref{L: Folland lemma}) and remembering that
$\delta_n$ is symmetric for all $n$,
\[
    (T\ast\delta_n)[\tau(x)]=
    T[(\delta_n\ast\tau)(x)]\to T[\tau(x)]
    \]
    by (Lemma~\ref{L: deltan*tau->tau}) and the
continuity of $T$. $\blacktriangle$

    \section{Embedding of Schwartz Distributions
     in $^*\mathcal{E}(\mathbb{R}^{d})$}\label{S:
Embedding of Schwartz Distributions
     in *E}
Finally, we are ready to define the embedding $\Sigma$ of
$\mathcal{D}'(\mathbb{R}^{d})$ into $^*
\mathcal{E}(\mathbb{R}^{d})$ as follows:
\[
    \Sigma(T) = \langle T \ast \varphi \rangle
    \]
    By (theorem~\ref{T: F*phi in E}),
    $\Sigma(T)\in\,^*\mathcal{E}(\mathbb{R}^d)$. From
    the definition of the convolution, it is clear
that
    $\Sigma$ is linear. It remains to prove that
$\Sigma$ is
    injective.

\begin{lemma}\label{L: embedding}
$\Sigma$ is injective.
\end{lemma}
\emph{Proof} Since $\Sigma$ is linear, it is enough to show that
$\Sigma(T) = 0$ implies $T=0$.\\
If $\Sigma(T) = 0$, we have that $T \ast \varphi = 0$ a.e. That
is, $\Phi = \{\varphi \, \vert \, T\ast \varphi = 0 \} \in
\mathcal{U}$. Thus $\varnothing \ne \Phi \cap \mathcal{B}_n \in
\mathcal{U}$ for each $n$, where $\mathcal{B}_n$ are the basic
sets. Therefore we can construct a sequence $\{\varphi_n\}$ such
that $\varphi_n \in \Phi \cap \mathcal{B}_n$ for each $n$. Then by
(Theorem~\ref{T: regularization}), we have that $T = 0$ since
$T\ast \varphi_n = 0$ for every $n$. $\blacktriangle$

\begin{theorem}\label{T: embedding}
\begin{itemize}
    \item[(i)] $\langle P \ast \varphi \rangle =\,
^*P$ for every
    polynomial $P \in \mathbb{C}[x_1, \ldots, x_d]$
    \item[(ii)] $\langle f \ast \varphi \rangle
    \overset{\rho}{=}\, ^*f$ for all
$f\in\mathcal{E}(\mathbb{R}^d)$
\end{itemize}
\end{theorem}
\emph{Proof} (i) Let $P \in \mathbb{C}[x_1, \ldots, x_d]$ be a
polynomial of degree $p$. By the Taylor formula,
\[
    P(x-t) = P(x) + \sum _{|\alpha | = 1}^{p} \frac
    {(-1^{|\alpha|})\partial^{\alpha}P(x)}{\alpha !}
t^{\alpha}
    \]
It follows that for every test function $\varphi$ and $x\in
\mathbb{R}^d$
\[
    (P \ast \varphi)(x) = \int P(x-t)\varphi(t) dt =
P(x) \int
    \varphi(t) dt  + \sum _{|\alpha | = 1} ^{p} \frac
    {(-1^{|\alpha|})\partial^{\alpha}P(x)}{\alpha !}
\int
    t^{\alpha} \varphi (t) dt
    \]
Notice that if $\varphi\in\mathcal{B}_n$ for some $n\geq p$, then
$\int \varphi (t) dt =1$ and $\int t^{\alpha} \varphi (t) dt = 0$,
$|\alpha| = 1, 2, \ldots, p$. Thus we have
\[
    \mathcal{B}_n \subseteq \{\varphi \, | \, P\ast
\varphi = P\}
    \]
implying that $P \ast \varphi = P$ a.e. as
required.\\\\
(ii) Let $\xi \in \mathcal{F}(^*\mathbb{R}^d)$, $n\in \mathbb{N}$,
and $\alpha$ be a multi-index. We have to show that $\vert
\partial ^{\alpha}(f\ast\varphi)(\xi) -
\partial^{\alpha}f(\xi)
\vert < \rho^{n} $. We will show this for the case $\alpha = 0$,
the general result will follow.\\
Since $st(\xi)\in\mathbb{R}^{d}$ we can find an open relatively
compact set $\mathcal{O}\subset\mathbb{R}^{d}$ such that
$st(\xi)\in\mathcal{O}$ and by (Robinson~\cite{aRob66}, p.90
Theorem 4.1.4) $\xi\in \, ^*\mathcal{O}$ and hence $\xi \in \,
^*\overline{\mathcal{O}}$.\\
As before, the Taylor formula gives

\[ f(x-t) = f(x) + \sum _{|\alpha | = 1}^{n} \frac
    {(-1^{|\alpha|})\partial^{\alpha}f(x)}{\alpha !}
t^{\alpha}
    + \sum_{|\alpha|=n+1} \frac
{(-1)^{|\alpha|}\partial^{\alpha}f(\eta(x,t))}{\alpha!}
    t^{\alpha}
    \]
where $\eta(x,t)$ is a point in $\mathbb{R}^d$ "between $x$ and
$t$". It follows that for every $\varphi \in
\mathcal{D}(\mathbb{R}^d)$
\[
    (f \ast \varphi)(x) = f(x)\int \varphi(t) dt +
    \sum_{|\alpha|=1}^{n} \frac
    {(-1)^{|\alpha|}\partial^{\alpha}f(x)}{\alpha !}
\int
    t^{\alpha} \varphi(t) dt +
    \]
\[
    \sum_{|\alpha|=n+1} \int \frac
    {(-1)^{|\alpha|}
\partial^{\alpha}f(\eta(x,t))}{\alpha !}
    t^{\alpha} \varphi(t) dt
    \]
Letting
\[
    M\overset{def}{=} 2 \sum_{|\alpha|=n+1} \sup_{x\in
K}
    \sup_{t\in\supp(\varphi)} \left|
    \frac{\partial^{\alpha}f(\eta(x,t))}{\alpha
!}\right|
    \]
we have that $M\rho^{n+1}<\rho^n$ since $M\in\mathbb{R}$ and
$\rho$ is a positive infinitesimal. In other words, if $\rho =
\langle R_\varphi \rangle $, then $M R_\varphi ^{n+1} < R_\varphi
^n$ a.e.\\
By the properties of $\mathcal{B}_n$, it follows that for a.e.
$\varphi$, $x\in K$,
\begin{align}\notag
    & |(f\ast \varphi)(x)-f(x)| < \sup_{x\in K} |(f
\ast
    \varphi)(x) - f(x)| \\ \notag
    & < \sum_{|\alpha|=n+1} \sup_{x\in K} \sup_{t \in
    \supp(\varphi)} \left|
    \frac{\partial^{\alpha}f(\eta(x,t))}{\alpha !}
\right| \left(
    \sup_{t\in\supp(\varphi)} \Vert t
    \Vert^{n-1}\right)\int|\varphi(t)|dt \\ \notag
    & \leq M R_{\varphi}^{n+1} < R_\varphi ^n
\end{align}
Finally, since $\xi \in \, ^*K$, we have that $\vert
 (f\ast\varphi)(\xi) - f(\xi)
\vert < \rho^{n} $, as required. The general result follows from
the case $\alpha = 0$ and the fact that $\partial^{\alpha} (f \ast
\varphi) = (\partial ^{\alpha} f) \ast \varphi$. $\blacktriangle$

\section{Conservation Laws in $^*\mathcal{E}(\Omega)$
and the Hopf Equation}\label{S: Integration in *E} The embedding
in the previous section is done deliberately, with the intent of
showing that $^*\mathcal{E}$ is a natural extension of
$\mathcal{D}'$ and an appropriate setting for the study of weak
solutions to non-linear partial differential equations, an
abundance of which arise from the conservation law of physics.

\begin{theorem}[Conservation Laws in
$^*\mathcal{E}(\Omega)$]\label{T: Conservation Laws} Let
$L\in\mathbb{R}_+\cup\{\infty\}$,
$F\in\mathcal{C}^\infty(\mathbb{C})$ and let $^*F$ be the
non-standard extension of $F$. Let $u\in{^*\mathcal{E}(\Omega)}$,
where $\Omega=(0, L)\times(0, \infty)$. Then the following are
equivalent:

{\bf (i)}\quad $u_t(x, t)+\left[^*F(u(x, t))\right]_x=0$ for all
$x, t\in{^*\mathbb{R}},\; 0<x<L,\; t>0$.

{\bf (ii)} \, $u_t(x, t)+{^*F^\prime\left(u(x, t)\right)}\, u_x=0$
for all $x, t\in{^*\mathbb{R}},\; 0<x<L,\; t>0$.

{\bf (iii)}\; $\frac{d}{dt}\int_a^b\, u(x, t)\, dx= {^*F(u(a,
t))}-{^*F(u(b, t))}$ for every $a, b, t\in{^*\mathbb{R}},\;
0<a<b<L,\; t>0$.
\end{theorem}
\begin{remark} The term ``conservation law'' is due to
(iii) which in a classical setting is given by
\[
\frac{d}{dt}\int_a^b\, u(x, t)\, dx= {F(u(a, t))}-{F(u(b, t))}.
\]
Here $u(x, t)$ stands for the density of a physical quantity (the
density of the mass of a fluid, the density of the heat energy,
etc.) in a rod of length $L$ and $F(u(x, t))$ stands for the flux
of the quantity from left to right through the $x$-cross section.
Then the above equality expresses the conservation of this
quantity in any $(a, b)$-segment of the rod. Recall that,
according to the classical theory, (i)-(iii) are equivalent for
solutions $u$ in the class $\mathcal{C}^2(\Omega)$ and for all $x,
t, a, b\in\mathbb{R}$ in the corresponding intervals. The proof
which follows can be generalized (without new complications) in
the case of more complicated flux $F(u, u_x)$ or even $F(u, u_x,
u_{xx})$.
\end{remark}
\emph{Proof} (i) $\Leftrightarrow$ (ii): The equivalency between
(i) and (ii) follows immediately from the fact that the partial
differentiation and extension mapping $*$ commute in
$^*\mathcal{E}(\Omega)$ and the fact that $^*\mathcal{E}(\Omega)$
is a differential algebra (with Leibniz rule for differentiation
of products and chain rule). So, we have
\[
\left[^*F(u(x, t))\right]_x={(^*F)^\prime(u(x, t))\,u_x(x,
t)}={^*F^\prime(u(x, t))\,u_x(x, t)},
\]
as required.

    (i) $\Rightarrow$ (iii): We have  $a=\left<
a_\varphi\right>, b=\left< b_\varphi\right>$ and $u=\left<
u_\varphi\right>$ for some families of real numbers $(a_\varphi),
(b_\varphi)\in\mathbb{R}^{\mathcal{D}(\mathbb{R}^2)}$ and some
family of smooth functions
$(u_\varphi)\in\mathcal{E}(\Omega)^{\mathcal{D}(\mathbb{R}^2)}$.
We have $\Phi=\{\varphi\mid a_\varphi< b_\varphi\}\in\mathcal{U}$
since $a<b$ in $^*\mathbb{R}$, by assumption. Thus (involving the
classical arguments in the framework of $\mathcal{E}(\Omega)$) we
have $\Phi\subseteq\Phi_1$, where
\[
\Phi_1=\{\varphi\mid \frac{d}{dt}\int_{a_\varphi}^{b_\varphi}\,
u_\varphi(x, t)\, dx= {F(u_\varphi(a_\varphi,
t))}-{F(u_\varphi(b_\varphi, t))} \text{\; for all\;}
t\in\mathbb{R}_+\}.
\]
The latter implies $\Phi_1\in\mathcal{U}$ which implies (iii), as
required, after transferring the result from representatives to
the corresponding equivalence classes.

    (i) $\Leftarrow$ (iii): Suppose (on the contrary)
that there exist $\xi, \tau\in{^*\mathbb{R}},\, 0< \xi< L,\,
\tau>0$, such that $u_t(\xi, \tau)+\left[^*F(u(\xi,
\tau))\right]_x\not=0$ in $^*\mathbb{C}$. We have
$\xi=\left<\xi_\varphi\right>$ and
$\tau=\left<\tau_\varphi\right>$ for some $(\xi_\varphi),
(\tau_\varphi)\in\mathbb{R}^{\mathcal{D}(\mathbb{R}^2)}$. We
denote $\Phi=\{\; \varphi\mid (u_{\varphi})_t(\xi_\varphi,
\tau_\varphi)+\left[F(u_\varphi(\xi_\varphi,
\tau_\varphi))\right]_x\not=0\; \}$ and observe that
$\Phi\in\mathcal{U}$ (by our assumption). Also, we let
\begin{align}\notag
\Phi_\varphi=\{\; (\alpha, \beta,\gamma)\in\mathbb{R}^3\mid &\;
\frac{d}{dt}\int_\alpha^\beta u_{\varphi}(x, \gamma)\, dx\not=
F(u_\varphi(\alpha, \gamma))- F(u_\varphi(\beta,\gamma)),\\\notag
&\; 0<\alpha<\beta<L,\; \gamma>0\;\},
\end{align}
and observe that $\Phi_\varphi\not=\varnothing$ for all
$\varphi\in\Phi$ (by the classical theory in the framework of
$\mathcal{E}(\Omega)$). By axiom of choice, there exist families
$(a_\varphi), (b_\varphi), (\gamma_\varphi)\in\mathbb{R}^{\Phi}$
such that $(a_\varphi, b_\varphi, \gamma_\varphi)\in\Phi_\varphi$
for all $\varphi\in\Phi$. If
$\varphi\in\mathcal{D}(\mathbb{R}^2)\setminus\Phi$, we define
$(a_\varphi, b_\varphi, \gamma_\varphi)$  anyhow (say, by
$a_\varphi=b_\varphi=\gamma_\varphi=1$). These families (of real
numbers) determine the non-standard real numbers $a=\left<
a_\varphi\right>, b=\left< b_\varphi\right>$ and
$t=\left<\gamma_\varphi\right>$. Next, we observe that
$\Phi\subseteq\Psi$ (by the definition of $\Phi_\varphi$), where
\begin{align}\notag
\Psi=\{\; \varphi\mid\; & \frac{d}{dt}\int_{a_\varphi}^{b_\varphi}
u_{\varphi}(x, \gamma_\varphi)\, dx\not= F(u_\varphi(a_\varphi,
\gamma_\varphi))- F(u_\varphi(b_\varphi,\gamma_\varphi)),\\\notag
                                & 0<a_\varphi<
b_\varphi< L,\; \gamma_\varphi>0  \}.
\end{align}
Next, $\Phi\in\mathcal{U}$ implies $\Psi\in\mathcal{U}$ which
implies
\[
\frac{d}{dt}\int_{a}^{b} u(x, t)\, dx\not= F(u(a,t))- F(u(b,t)),\;
0<a< b< L,\; t>0,
\]
in the framework of $^*\mathbb{C}$, contradicting (iii). The proof
is complete. $\blacktriangle$

\begin{example}[Hopf Equation] To appreciate the
result of Theorem~\ref{T: Conservation Laws} we recall that
(i)-(iii) might or might not be equivalent in classes of classical
functions and Schwartz distributions larger than
$\mathcal{C}^2(\Omega)$, where the (important for the theory and
applications) {\bf fundamental solutions} and {\bf shock wave
solutions} belong. For a discussion we refer to (J. David
Logan~\cite{Logan}, p. 309-310). Here is an example: Let
$F(u)=\frac{1}{2}u^2$ and $L=\infty$, so we have
$\Omega=\mathbb{R}^2_+$. In this case (i)-(iii) become:

{\bf (i)}\quad $u_t(x, t)+[\frac{1}{2}u(x, t)^2]_x=0$ for all $x,
t\in{\mathbb{R}},\; x>0,\; t>0$.

{\bf (ii)} \, $u_t(x, t)+u(x, t)u_x(x, t)=0$ ({\bf Hopf equation})
for all $x, t\in{\mathbb{R}},\; x>0,\; t>0$.

{\bf (iii)}\; $\frac{d}{dt}\int_a^b\, u(x, t)\, dx=
\frac{1}{2}\left[u^2(a, t)-u^2(b, t)\right]$ for every $a, b,
t\in{\mathbb{R}},\; 0<a<b,\; t>0$,

\noindent respectively. Let $v\in\mathbb{R}_+$, $H$ be the
Heaviside step function and let $u(x, t)=2vH(x-vt)$ be a {\bf
shock wave}. The next analysis shows that (i), (ii) and (iii) are
not equivalent in the spaces of classical functions and Schwartz
distributions:

    {\rm (i)} Since
$u=2vH(x-vt)\notin\mathcal{C}^2(\mathbb{R}^2_+)$ this function can
not be a classical solution of (i). However, $u=2vH(x-vt)$ is a
(generalized) solution of (i) in the framework of the class of
Schwartz distributions $\mathcal{D}^\prime(\mathbb{R}^2_+)$.
Indeed, for the first term of (i) we have $u_t(x, t)=
-2v^2\delta(x-vt)$, where $\delta(x)$ is the Dirac delta function.
For the second term we have $[\frac{1}{2}u(x,
t)^2]_x=[\frac{4v^2}{2}H(x-vt)^2]_x=[2v^2H(x-vt)]_x=2v^2\delta(x-vt)$.
Thus $u(x, t)=2vH(x-vt)$ is a (generalized) solution of (i).  We
should notice that $u(x, t)=2vH(x-vt)$  is also a weak solution of
(i) in the framework of $\mathcal{L}_{loc}(\mathbb{R}^2_+)$ (see
the remark below).

    {\rm (ii)} $u(x, t)=2vH(x-vt)$ is clearly not a
solution of (ii) in classical sense. Neither is it a (generalized)
solution of (ii) in the class $\mathcal{D}^\prime(\mathbb{R}^2_+)$
because the term $uu_x=4v^2H(x-vt)\delta(x-vt)$ does not make
sense within $\mathcal{D}^\prime(\mathbb{R}^2_+)$ (recall that
there is no multiplication in
$\mathcal{D}^\prime(\mathbb{R}^2_+)$).

    {\rm (iii)} For the LHS of (iii) we have
$\frac{d}{dt}\int_a^b\, u(x, t)\,
dx=\frac{d}{dt}\int_a^b\,2vH(x-vt)\,
dx=2v\frac{d}{dt}\int_{a-vt}^{b-vt}\,H(x)\, dx =-2v^2H(vt-a)$. For
the RHS of (iii) we have \newline $\frac{1}{2}\left[u^2(a,
t)-u^2(b,
t)\right]=2v^2\left[H(a-vt)-H(b-vt)\right]=-2v^2H(vt-a)$. Thus
$u(x, t)=2vH(x-vt)$ is a solution of (iii).

\end{example}

\begin{remark}[Weak Solution] Suppose the
$u\in\mathcal{L}_{loc}(\Omega)$ is a solution of
$u_t(x,t)+\left[F(u(x, t))\right]_x=0$ in the framework of
$\mathcal{D}^\prime(\Omega)$ (that means that both $u_t(x, t)$ and
$\left[F(u(x, t))\right]_x$ are in $\mathcal{D}^\prime(\Omega)$).
These solutions are often called {\bf weak solutions} because they
satisfy the {\bf weak equality}:
\[
\iint_{\Omega}\left[u(x,t)\tau_t(x, t)+F(u(x, t))\tau_x(x,
t)\right]\,dx\, dt=0,
\]
for all test functions $\tau\in\mathcal{D}(\Omega)$. In
$\mathcal{D}^\prime(\Omega)$ we have:
\[
\iint_{\Omega}\left[u(x,t)\tau_t(x, t)+F(u(x, t))\tau_x(x,
t)\right]\,dx\, dt=\left<u_t(x,t)+\left[F(u(x, t))\right]_x,
\tau(x, t)\right>,
\]
where $\left<\; ,\; \right>$ stands for the pairing between
$\mathcal{D}^\prime(\Omega)$ and $\mathcal{D}(\Omega)$. Thus every
(generalized) solution in $\mathcal{D}^\prime(\Omega)$ is also
weak solution in $\mathcal{L}_{loc}(\Omega)$. In particular the
function $u(x, t)=2vH(x-vt)$ in the example above is both a
generalized and a weak solution of (i).
\end{remark}

\section{Generalized Delta-like Solution of the Hopf
Equation}\label{S: Generalized Delta-Like Solution of Hopf
Equation}

In this section we will prove the existence of a delta-like weak
solution to the Hopf equation of the type
$\overset{\rho}{\simeq}$. That is, we are looking for a function
of the form
\[
    u(x,t) = u_0 + \frac{A}{\rho}
\Theta\left(\frac{x-vt}{\rho}\right)
    \]
    where $\Theta \in \, ^*\mathcal{S}(\mathbb{R})$
($\mathcal{S}(\mathbb{R})$
    is the
    class of rapidly decreasing functions, such as
$e^{-x^2}$, $^*\mathcal{S}$ is
    its non-standard extension, defined similarly to
$^*\mathcal{E}$),
    $\int \Theta (x) dx = 1$, $u_0, A, v\in
\mathcal{M}_\rho(^*\mathbb{R})$
     (we consider A to be
    the amplitude of the soliton and $v$ its
velocity), and for all $t>0$,
\[
    u_t + uu_x \overset{\rho}{\simeq} 0
    \]
    That is, for all $t>0$,
$\tau\in\mathcal{D}'(\mathbb{R})$,
\[
    \int [u_t + uu_x] \tau(x)dx \overset{\rho}{=} 0
    \]
In addition, we would like this function to satisfy the
conservation law, so that for all $a, b\in\mathbb{R}$, $t>0$,
\[
    \frac{d}{dt} \int_a ^b u(x,t) dx \overset{\rho}{=}
    \frac{1}{2}[u^2(a,t)-u^2(b,t)]
    \]
Calculating $u_t + uu_x$ we get
\[
    u_t + uu_x = -\frac{Av}{\rho^2}\Theta '
    \left(\frac{x-vt}{\rho}\right) + \frac{u_0
A}{\rho^2}\Theta '
    \left(\frac{x-vt}{\rho}\right) +
\frac{A^2}{\rho^3}\Theta\left(\frac{x-vt}{\rho}\right)\Theta '
    \left(\frac{x-vt}{\rho}\right)
    \]
Simplifying and letting
$\Theta\left(\frac{x-vt}{\rho}\right)\Theta '
\left(\frac{x-vt}{\rho}\right)=\frac{\rho}{2}\left(\Theta^2
    \left(\frac{x-vt}{\rho}\right)\right)_x$ gives us

\begin{align}\notag
    & \int [u_t + uu_x] \tau(x) dx \\ \notag
    & = \frac{(u_0 - v)A}{\rho^2} \int \Theta
    '\left(\frac{x-vt}{\rho}\right)\tau(x) dx +
\frac{A^2}{2\rho^2}\int\left(\Theta^2\left(\frac{x-vt}{\rho}\right)\right)_x
     \tau(x) dx
\end{align}
Integrating by parts and making the substitution $y =
\frac{x-vt}{\rho}$ gives

\[
   = \int \left[(v-u_0)A\Theta(y) -
\frac{A^2}{2\rho}\Theta^2(y)\right]\tau ' (vt+\rho y) dy
    \]
Finally, using the Taylor formula for $\tau'(vt+\rho y)$, we have
that for each $m\in\mathbb{N}$,
\[
    = \sum_{n=0}^m \int \left[(v-u_0)A\Theta(y) -
    \frac{A^2}{2\rho}\Theta^2(y)\right]y^n
        \frac{\tau^{(n+1)}(vt)}{n!}\rho^n \, dy +
R_m(\tau)
\]
where the remainder term is
\[
    R_m(\tau) =  \rho^{m+1}\int \left[(v-u_0)A
\Theta(y)-
    \frac{A^2}{2\rho}\Theta^2(y)\right]
    \frac{\tau^{(m+2)}(\eta(y,t))}{(m+1)!}y^{m+1}dy
    \]
We would like to find a function $\Theta$ such that for every $m$,
\[
    \int\left[(v-u_0)A\Theta(y) -
    \frac{A^2}{2\rho}\Theta^2(y)\right]y^n dy = 0, \,
0\leq n\leq m
    \]
and $\vert R_m(\tau) \vert < \rho^{m+k}$ (for some
fixed k).\\
When $m=0$, we have that
\[
    A = \frac{2\rho(v-u_0)}{\int\Theta^2(y) dy}
    \]
    (remembering that $\int \Theta(y) dy = 1$)
Replacing this value of A, we have that for every $m$,
\[
    \int \Theta^2(y) dy \int \Theta (y) y^n dy = \int
\Theta^2(y)
    y^n dy, \, 0\leq n \leq m
    \]
Define
\[
    S_m = \{ f\in \mathcal{S} \, : \, \int f(x)x^n dx
=
    \frac{\int f^2(x)x^n dx}{\int f^2(x) dx} , \,
0\leq n\leq m\}
    \]
For each $m$, $S_m$ is non-empty by (M.
Radyna~\cite{mRadynaICGF2000} p. 275).\\\\
Now let
\begin{align}\notag
    \overline{S}_m =  \{ f\in\, ^*S_m \, : & f(0)=0
\\ \notag
                                          &
|\ln\rho|^{-1}\int |f(x)x^n| < 1 \\ \notag
                                          &
|\ln\rho|^{-1}\int
|f^2(x)x^n| dx < 1, \, 0\leq n \leq m \}\\
                     \notag
                     \end{align}
The standard functions in $^*S_m$ certainly satisfy the second and
third conditions, since their integrals will be standard (and
therefore finite) and $|\ln\rho|^{-1}$ is infinitely small. As for
the first condition, we can say with certainty that if a function
 $f\in\,^*S_m$ has at least one zero, say $f(-k)=0$,
then
 $g(x)\overset{def}{=}f(x-k)\in\,^*S_m$ and $g(0)=0$
by the
 following lemma:

 \begin{lemma}\label{L: invariance under translation}
 Suppose $f(x)$ satisfies
 \[
    \int f^2(x) dx \int f (x) x^n dx = \int f^2(x)
    x^n dx
    \]
    Then $g(x) = f(x-k)$ also satisfies
 \[
    \int g^2(x) dx \int g (x) x^n dx = \int g^2(x)
    x^n dx
    \]
 \end{lemma}
\emph{Proof} Substituting $y=x-k$, we get
\begin{align} \notag
    \int g^2(x) dx \int g (x) x^n dx &= \int f^2(y)dy
\int
    f(y)(y+k)^n dy \\ \notag
                                     &=\sum_{j=0}^n
                                     \begin{pmatrix}
                                     n \\
                                     j
                                     \end{pmatrix}
                                     k^{n-j} \int
f^2(y)dy \int
                                     f(y) y^j dy \\
\notag
                                     &=\sum_{j=0}^n
                                     \begin{pmatrix}
                                     n \\
                                     j
                                     \end{pmatrix}
                                     k^{n-j}
                                     \int f^2(y) y^j
dy \\ \notag
                                     &=\int f^2(y)
\sum_{j=0}^n
                                     \begin{pmatrix}
                                     n \\
                                     j
                                     \end{pmatrix}
                                     y^jk^{n-j} dy \\
\notag
                                     &=\int f^2(y)
(y+k)^ndy \\
                                     \notag
                                     &=\int f^2(x-k)
x^n dx\\ \notag
                                     &=\int g^2(x) x^n
dx
                                     \quad \quad
\blacktriangle
                                     \end{align}
Thus, if at least one function in $^*S_m$ has at least one zero,
then $\overline{S}_m$ will be non-empty. In addition, the sets
$\overline{S}_m$ are internal and $\overline{S}_0 \supset \,
\overline{S}_1 \supset \ldots$. Therefore, by the saturation
principle, there exists a function
$\Theta(x)\in\bigcap_{n=0}^\infty \overline{S}_n$.\\
Noting that $\int |\Theta(x)x^n| dx$ and $\int |\Theta^2(x)x^n|dx$
are at most $\mathcal{C}_\rho(\mathbb{^*C})$, we have that for
this $\Theta$,
\[
    |R_m(\tau)| <
\frac{\rho^{m+1}}{(m+1)!}\sup_{x\in\mathbb{R}}|\tau^{(m+2)}(x)|
    \int\left|
\left[(v-u_0)A\Theta(y)-\frac{A^2}{2\rho}\Theta^2(y)\right]y^n\right|
    dy < \rho^{m+k}
    \]
    where $k$ is some real constant that depends on
$v, u_0,$ and $A$.
    Therefore
$u(x,t)=u_0+\frac{A}{\rho}\Theta\left(\frac{x-vt}{\rho}\right)$
     satisfies the Hopf equation weakly, in the
    sense of  $\overset{\rho}{\simeq}$.\\
    If, in addition, it is true that $\Theta(0)=0$,
then $u(x,t)$
    will also satisfy the conservation law:\\
    for all $a, b\in\mathbb{R}$, $t>0$,
\[
    \frac{d}{dt} \int_a ^b u(x,t) dx \overset{\rho}{=}
    \frac{1}{2}[u^2(a,t)-u^2(b,t)]
    \]
Let us prove this by first calculating the left side:
\begin{align}\notag
    \frac{d}{dt} \int_a^b \left[u_0 +
\frac{A}{\rho}\Theta\left(\frac{x-vt}{\rho}\right)\right]dx &=
    \frac{d}{dt}\left[u_0(b-a)+\frac{A}{\rho}
\int_a^b\Theta\left(\frac{x-vt}{\rho}\right)dx\right]
\\
    \notag
&=\frac{d}{dt}\left[A\int_{\frac{a-vt}{\rho}}^{\frac{b-vt}{\rho}}\Theta(y)dy\right]
    \\\notag
&=A\left[\Theta\left(\frac{b-vt}{\rho}\right)\left(-\frac{v}{\rho}\right)-
\Theta\left(\frac{a-vt}{\rho}\right)\left(-\frac{v}{\rho}\right)\right]\\\notag
&=A\frac{v}{\rho}\left[\Theta\left(\frac{a-vt}{\rho}\right)
    -\Theta\left(\frac{b-vt}{\rho}\right)\right]
    \end{align}
Since $\Theta\in\,^*\mathcal{S}$, if $a\neq vt$ and $b\neq vt$
then this quantity vanishes. However, if $a=vt$ or $b=vt$ (but not
both) then we have that the left side equals
$\pm\frac{Av}{\rho}\Theta(0)$, respectively. If $\Theta(0)=0$,
this also vanishes.\\
Calculating the right side, we have:
\begin{align}\notag
    &\frac{1}{2} \left[\left[u_0 +
\frac{A}{\rho}\Theta\left(\frac{a-vt}{\rho}\right)\right]^2
    -\left[u_0 +
\frac{A}{\rho}\Theta\left(\frac{b-vt}{\rho}\right)\right]^2\right]\\\notag
    &=\frac{1}{2} u_0^2+
\frac{u_0A}{\rho}\Theta\left(\frac{a-vt}{\rho}\right)
+\frac{A^2}{2\rho^2}\Theta^2\left(\frac{a-vt}{\rho}\right)\\\notag
&-\frac{1}{2}u_0^2-\frac{u_0A}{\rho}\Theta\left(\frac{b-vt}{\rho}\right)
-\frac{A^2}{2\rho^2}\Theta^2\left(\frac{b-vt}{\rho}\right)
\\\notag
&=\frac{u_0A}{\rho}\left[\Theta\left(\frac{a-vt}{\rho}\right)
        -\Theta\left(\frac{b-vt}{\rho}\right)\right]
        +
\frac{A^2}{2\rho^2}\left[\Theta^2\left(\frac{a-vt}{\rho}\right)
        -\Theta^2\left(\frac{b-vt}{\rho}\right)\right]
\end{align}
Again, we have that if $a\neq vt$ and $b\neq vt$ then this
quantity vanishes. If $a=vt$ or $b=vt$ (but not both) then the
right side equals
$\pm\left(\frac{u_0A}{\rho}\Theta(0)+\frac{A^2}{2\rho^2}\Theta^2(0)\right)$,
respectively. Here also, if $\Theta(0)=0$ the right side vanishes,
and so the conservation law holds.\\\\
In conclusion, we may make some conjectures based on the relation
\[
    A = \frac{2\rho(v-u_0)}{\int\Theta^2(y) dy}
    \]
There are many possibilities here, but if we assume (for
simplicity) that $\int\Theta^2(y) dy$ is finite (and not
    infinitesimal) and $u_0=0$, then there are at
least the following two particular cases:
        \begin{itemize}
            \item[(i)] u has infinitesimal amplitude
with finite or
            infinitely large velocity, resembling a
small signal,
            or
            \item[(ii)] u has non-infinitesimal,
finitely large amplitude,
            and infinitely large velocity, resembling
an
            explosion.
            \end{itemize}
\begin{remark}[Connection with Perturbation Theory]
The closest to our result is the work by M.
Radyna~\cite{mRadynaICGF2000} in the framework of perturbation
theory. M. Radyna proves the following result: For every
$n\in\mathbb{N}$ there exists a function
$\Theta_n\in\mathcal{S}(\mathbb{R})$ such that the function $u(x,
t)=\frac{A}{\epsilon}\Theta_n(\frac{x-vt}{\epsilon})$ satisfies:
\[
\left|\int_{\mathbb{R}}[u_t(x, t)+u(x, t)u_x(x, t)]\tau(x)\, dx
\right|<\epsilon^n,
\]
for every test function $\tau\in\mathcal{D}(\mathbb{R})$,  every
$t\in\mathbb{R}$ and all sufficiently small
$\epsilon\in\mathbb{R}$. For comparison we mention the following:

    {\bf (a)} Instead of a small real parameter
$\epsilon$ we use a proper positive infinitesimal $\rho$. Our
framework is, of course, quite different from M. Radyna's theory.

    {\bf (b)} In contrast to M. Radyna's result, we
have proved the existence of a function
$\Theta\in{^*\mathcal{S}(\mathbb{R})}$ (not depending on $n$) such
that the function $u(x,
t)=\frac{A}{\rho}\Theta(\frac{x-vt}{\rho})$ satisfies:
\[
\left|\int_{\mathbb{R}}[u_t(x, t)+u(x, t)u_x(x, t))\tau(x)\, dx
\right|<\rho^n,
\]
for every test function $\tau\in\mathcal{D}(\mathbb{R})$,  every
$t\in\mathbb{R}$ and {\bf for all} $n\in\mathbb{N}$.

\end{remark}
\newpage


\end{document}